\documentclass{gtart}
\usepackage{amssymb, amsmath, psfrag, epsfig, graphs}

\newtheorem{thm}{Theorem}[section]
\newtheorem{cor}[thm]{Corollary}
\newtheorem{lem}[thm]{Lemma}
\newtheorem{prop}[thm]{Proposition}
\newtheorem{defn}[thm]{Definition}

\theoremstyle{remark}
\newtheorem{rem}[thm]{Remark}
\numberwithin{equation}{section}

\newcommand{\de}{\delta}

\newcommand{\Si}{\Sigma}

\newcommand{\bD}{\mathbf D}

\newcommand{\RR}{\mathcal R}

\newcommand{\N}{\mathbb N}
\newcommand{\Z}{\mathbb Z}
\newcommand{\Q}{\mathbb Q}

\newcommand{\del}{\partial}

\newcommand{\hra}{\hookrightarrow}

\DeclareMathOperator{\Aut}{Aut}

\begin{document}

\title{Lens spaces, rational balls and the ribbon conjecture} 

\author{Paolo Lisca}
\address{Dipartimento di Matematica ``L. Tonelli''\\ 
Largo Bruno Pontecorvo, 5\\
Universit\`a di Pisa \\
I-56127 Pisa, ITALY} 


\begin{abstract}
We apply Donaldson's theorem on the intersection forms of definite
4--manifolds to characterize the lens spaces which smoothly bound
rational homology 4--dimensional balls. Our result implies, in
particular, that every smoothly slice 2--bridge knot is ribbon,
proving the ribbon conjecture for 2--bridge knots.
\end{abstract}

\maketitle

\section{Introduction}
\label{s:intro}

It is a well--known fact that every ribbon knot is smoothly slice. The
{\rm ribbon conjecture} states that, conversely, a smoothly slice knot
is ribbon. In this paper we prove that the ribbon conjecture holds for
2--bridge knots, deducing this result from a characterization of the
3--dimensional lens spaces which smoothly bound rational homology
4--dimensional balls (Theorem~\ref{t:main} below).

\begin{figure}[ht]
\begin{center}
\psfrag{-c1}{${\scriptstyle -c_1}$} 
\psfrag{c2}{${\scriptstyle c_2}$}
\psfrag{-c3}{${\scriptstyle -c_3}$} 
\psfrag{-cn-2}{${\scriptstyle -c_{n-2}}$} 
\psfrag{cn-1}{${\scriptstyle c_{n-1}}$}
\psfrag{-cn-1}{${\scriptstyle -c_{n-1}}$}
\psfrag{cn}{${\scriptstyle c_n}$} 
\psfrag{-cn}{${\scriptstyle -c_n}$} 
\psfrag{n even}{($n$ even)}
\psfrag{n odd}{($n$ odd)}
\psfrag{s}{\small $s$}
\psfrag{crossings}{$|s|$ crossings} 
\psfrag{s>0}{if $s\geq 0$}
\psfrag{s<0}{if $s\leq 0$} 
\includegraphics[height=8cm]{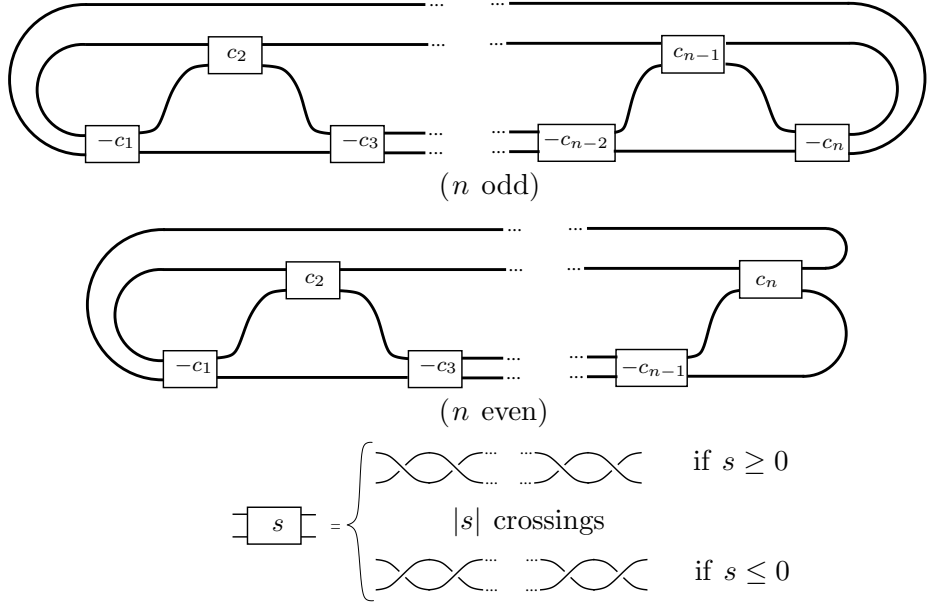}
\end{center}
\caption{The $2$--bridge link $L(c_1,\ldots,c_n)$}
\label{f:fig1}
\end{figure}

A link in $S^3$ is called~\emph{$2$--bridge} if it can be isotoped
until it has exactly two local maxima with respect to a standard
height function. Figure~\ref{f:fig1} represents the $2$--bridge link
$L(c_1,\ldots,c_n)$, where $c_i\in\Z$, $i=1,\ldots,n$.
Given coprime integers $p>q>0$ with    
\[
\frac pq = c_1 + \cfrac{1}{c_2 +
       \cfrac{1}{\ddots +
        \cfrac{1}{c_n}
}},
\quad c_i>0\quad\text{for $i=1,\ldots,n$},
\]
the 2--bridge link $K(p,q)$ is, by definition, $L(c_1,\ldots,c_n)$.
When $p$ is even, $K(p,q)$ is a 2--component link, when $p$ is odd
$K(p,q)$ is a knot. It is well--known~\cite[Chapter~12]{BZ} that
$K(p,q)$ and $K(p',q')$ are isotopic if and only if $p=p'$ and either
$q = q'$ or $qq'\equiv 1\pmod p$, and that every
$2$--bridge link is isotopic to some $K(p,q)$. Moreover, $K(p,p-q)$ is
isotopic to the mirror image of $K(p,q)$.

Now we recall what is known about 2--bridge knots with regard
to the ribbon conjecture. In order to do that, the following
definition is needed.
\begin{defn}\label{d:R}
Let $\Q_{>0}$ denote the set of positive rational numbers, and define
maps $f,g\co\Q_{>0}\to\Q_{>0}$ by setting, for $\frac pq\in\Q_{>0}$,
$p>q>0$, $(p,q)=1$,
\[
f\left(\frac pq\right) = \frac p{p-q},\qquad 
g\left(\frac pq\right )= \frac p{q'},
\]
where $p>q'>0$ and $qq'\equiv 1\pmod p$. Define $\RR\subset\Q_{>0}$ to
be the smallest subset of $\Q_{>0}$ such that $f(\RR)\subseteq\RR$,
$g(\RR)\subseteq\RR$ and $\RR$ contains the set of rational numbers
$\frac pq$ such that $p>q>0$, $(p,q)=1$, $p=m^2$ for some
$m\in\N$ and $q$ is of one of the following types:
\begin{enumerate}
\item
$mk\pm 1$ with $m>k>0$ and $(m,k)=1$;
\item
$d(m\pm 1)$, where $d>1$ divides $2m\mp 1$;
\item
$d(m\pm 1)$, where $d>1$ is odd and divides $m\pm 1$.
\end{enumerate}
\end{defn}

According to~\cite{Si}, Casson, Gordon and Conway showed that every
knot of the form $K(p,q)$ with $\frac pq\in\RR$ is ribbon. The
interior of any ribbon disk can be radially pushed inside the 4--ball
$B^4$ to obtain a smoothly embedded disk, and the 2--fold cover of
$B^4$ branched along a slicing disk for $K(p,q)$ is a smooth rational
homology ball with boundary the lens space $L(p,q)$. Therefore if
$K(p,q)$ is a knot (i.e.~if $p$ is odd) we have the following
implications:
\[
\frac pq\in\RR\ \Rightarrow\ K(p,q)\ \text{ribbon}\ 
\Rightarrow K(p,q)\ \text{smoothly slice}\ \Rightarrow
\  L(p,q) = \del W,
\]
where $W$ is a smooth 4--manifold with $H_*(W;\Q)\cong H_*(B^4;\Q)$.
Casson and Gordon~\cite{CG} observed that if $K(p,q)$ is a smoothly
slice knot then $p$ is a perfect square. Moreover, they proved that if
the 2--bridge knot $K(m^2,q)$ is ribbon then
\begin{equation}\label{e:CG}
\frac 2{m^2} \sum_{s=1}^{m^2-1} \cot(\frac{\pi s}{m^2}) 
\cot(\frac{\pi qs}{m^2})\sin^2(\frac{\pi rs}m) = \pm 1, 
\quad r=1,\ldots,m-1.
\end{equation}
Casson and Gordon~\cite[page 188]{CG} used Equations~\eqref{e:CG} to
show that if a $2$--bridge knot $K(m^2,q)$ is ribbon and $m\leq 105$
then $\frac{m^2}q$ belongs to $\RR$. Fintushel and
Stern~\cite[Theorem~6.1]{FS} proved that Equations~\eqref{e:CG} hold
under the assumption that $L(m^2,q)$ is the boundary of a smooth
rational homology ball $W$ with $H_2(W;\Z)$ without
2--torsion. In~\cite{OS} Owens and Strle used a result by Oszv\'ath
and Szab\'o~\cite[Theorem~9.6]{OSz} to find apriori different
obstructions for $K(m^2,q)$ to be smoothly slice, and verified that
for $m\leq 105$ these new obstructions give the same constraints as
Equations~\eqref{e:CG}. It is not known whether Equations~\eqref{e:CG}
imply that the knot $K(m^2,q)$ is smoothly slice.

The following is our main result.

\begin{thm}\label{t:main}
Let $p>q>0$ be coprime integers. Then, the following statements are
equivalent:
\begin{enumerate}
\item
The lens space $L(p,q)$ smoothly bounds a rational homology ball.
\item
There exist: 
\begin{enumerate}
\item
A surface with boundary $\Si$, homeomorphic to a disk if $p$ is odd and
to the disjoint union of a disk and a M\"obius band if $p$ is even;
\item
A ribbon immersion $i\co \Si\looparrowright S^3$ with 
$i(\del \Si)=K(p,q)$. 
\end{enumerate}
\item
$\frac pq$ belongs to $\RR$.
\end{enumerate}
\end{thm}

Theorem~\ref{t:main} immediately implies the following result, which 
settles the ribbon conjecture for 2--bridge knots.

\begin{cor}\label{c:main}
Let $p>q>0$ be coprime integers with $p$ odd. Then, the following
statements are equivalent:
\begin{enumerate}
\item
$\frac pq$ belongs to $\RR$
\item
$K(p,q)$ is a ribbon knot
\item 
$K(p,q)$ is a smoothly slice knot
\item
$L(p,q)$ smoothly bounds a rational homology ball.
\end{enumerate}
In particular, the ribbon conjecture holds for 2--bridge knots.
\qed\end{cor}

The proof of Theorem~\ref{t:main} is based on the following simple
idea. If a lens space $L(p,q)$ smoothly bounds a rational homology
ball $W(p,q)$, one can form a smooth negative definite 4--manifold
$X(p,q)$ by taking the union of $-W(p,q)$ with a canonical
4--dimensional plumbing $P(p,q)$ bounding $L(p,q)$. Since $X(p,q)$ is
negative definite, Donaldson's celebrated theorem~\cite{Do} implies
that the intersection lattice $Q_{X(p,q)}$ of $X(p,q)$ is isomorphic
to the standard diagonal intersection lattice $\bD^n$, where
$n=b_2(X(p,q))$. Therefore there is an embedding of intersection
lattices $Q_{P(p,q)}\hra\bD^n$, and since $-L(p,q)=L(p,p-q)$ smoothly
bounds the rational homology ball $-W(p,q)$, for some $n'$ there is an
embedding $Q_{P(p,p-q)}\hra \bD^{n'}$ as well. The existence of both
embeddings (it is easy to see that a single embedding is not enough)
gives constraints on the pair $(p,q)$ which eventually lead to the
proof of Theorem~\ref{t:main}. In spite of the simplicity of this
idea, the algebro--combinatorial machinery we must set up to work out
the above constraints is fairly complex and occupies
Sections~\ref{s:prelim}--\ref{s:strings} of the paper. Here is the
gist of what we do. We can write
\[
\frac pq = a_1 - \cfrac{1}{a_2 -
       \cfrac{1}{\ddots -
        \cfrac{1}{a_n}
}},
\qquad
\frac p{p-q} = b_1 - \cfrac{1}{b_2 -
       \cfrac{1}{\ddots -
        \cfrac{1}{b_{n'}}
}},
\]
for some integers $a_i,b_j\geq 2$ for $i=1,\ldots,n$,
$j=1,\ldots,n'$. It turns out (see Lemma~\ref{l:negsum}) that
\[
\sum_{i=1}^n (a_i-3) + \sum_{j=1}^{n'} (b_j-3) = -2,
\]
therefore, up to replacing $(p,q)$ with $(p,p-q)$, we may assume
\begin{equation}\label{e:condition}
\sum_{i=1}^n (a_i-3)<0. 
\end{equation}
After choosing a suitable set of generators of $H_2(P(p,q);\Z)$, the
embedding $Q_{P(p,q)}\hra\bD^n$ gives rise to a subset
$S=\{v_1,\ldots,v_n\}\subset\bD^n$ with
\[
v_i\cdot v_j =
\begin{cases}
-a_i\quad\text{if $|i-j|=0$},\\ 
1\quad\text{if $|i-j|=1$},\\ 
0\quad\text{if $|i-j|>1$}.
\end{cases}
\]
We call such subsets \emph{standard}. In
Sections~\ref{s:prelim}--\ref{s:strings} we study the standard subsets
of $\bD^n$ satisfying Equation~\eqref{e:condition}. In
Section~\ref{s:strings} we show that the string of integers
$(a_1,\ldots,a_n)$ associated to such a subset must belong (for a
fixed $n$) to a finite list which we describe explicitly. This gives
the constraints mentioned above. In Section~\ref{s:ribbon} we prove
the existence of ribbon surfaces for all the links required by
Theorem~\ref{t:main}~\footnote{The results of Section~\ref{s:ribbon}
were known previously for knots~\cite{Si} (although even in the case
of knots we were unable to recover all of them from~\cite{Si}). In
Section~\ref{s:ribbon} we give a self--contained account valid for
links and adapted to our conventions.}, and in Section~\ref{s:final}
we prove Theorem~\ref{t:main} using all the results obtained in the
previous sections. Each section starts with a brief outline
summarizing its purpose, contents and relationships with the other
sections.

{\bf Acknowledgments}: The author is grateful to Andrew Casson for
generous help, to Cameron Gordon for informative e--mail
correspondence and to the anonymous referee for useful comments. 

\section{First definitions and preliminary results}
\label{s:prelim}

\textbf{Outline.}  In this section we introduce definitions which will be
used throughout the paper. In particular, the concept of~\emph{good
subset} (see Definition~\ref{d:good}) is crucial in
Sections~\ref{s:p1>0},~\ref{s:p2>0} and~\ref{s:general}, while the
invariant $I(S)$ (see Definition~\ref{d:invariant}) is the key
quantity on which the proof of Theorem~\ref{t:main} is based. We also
prove Lemma~\ref{l:n=3}, which is the basis of the inductive process
used in the subsequent sections, and Lemma~\ref{l:negsum}, which will
be directly quoted in the proof of Theorem~\ref{t:main} in
Section~\ref{s:final}.

Let $\bD$ denote the intersection lattice $(\Z,(-1))$, and 
let $\bD^n$ be the orthogonal direct sum of $n$ copies of $\bD$. 
Fix generators $e_1,\ldots,e_n\in\bD^n$ such that
\[
e_i\cdot e_j=-\de_{ij},\quad i,j=1,\ldots,n. 
\]

Observe that the group of automorphisms $\Aut(\bD^n)$ contains the
reflections across each hyperplane orthogonal to an $e_i$ as well as
the all the trasformations determined by the permutations of
$\{e_1,\ldots,e_n\}$. Given a subset
$S=\{v_1,\ldots,v_n\}\subseteq\bD^n$, we define
\[
E^S_i:=\{j\in\{1,\ldots,n\}\ |\ v_j\cdot e_i\neq 0\},\quad
i=1,\ldots,n,
\]
\[
V_i:=\{j\in\{1,\ldots, n\}\ |\ e_j\cdot v_i\neq 0\},\quad
i=1,\ldots,n,
\]
and 
\[
p_i(S):=|\{j\in\{1,\ldots,n\}\ |\ |E^S_j|=i\}|,\quad
i=1,\ldots,n.
\]

Let $v_1,\ldots, v_n\in\bD^n$ be elements such that, for
$i,j\in\{1,\ldots,n\}$,
\begin{equation}\label{e:conds}
v_i\cdot v_j = 
\begin{cases}
-a_i\leq -2\quad& \text{if}\ i=j,\\
0\ \text{or}\ 1 \quad& \text{if}\ |i-j|=1,\\
0\quad& \text{if}\ |i-j|>1.
\end{cases}
\end{equation}
for some integers $a_i$, $i=1,\ldots,n$.

\begin{rem}\label{r:indep}
Elements $v_1,\ldots, v_n\in\bD^n$ satisfying
Conditions~\eqref{e:conds} are linearly independent over $\Z$. In
fact, it is easy to check that the associated intersection matrix
\[
Q := (v_i\cdot v_j)
\]
is nonsingular. The independence of $v_1,\ldots, v_n$ follows
immediately from the fact that
\[
Q= -MM^t, 
\]
where $M:=(m_{ij})$ is defined by $v_i=\sum_j m_{ij} e_j$.
\end{rem}

Let $S=\{v_1,\ldots, v_n\}\subseteq\bD^n$ be a subset which
satisfies~\eqref{e:conds}. We define the~\emph{intersection graph} of
$S$ as the graph having as vertices $v_1,\ldots,v_n$, and an edge
between $v_i$ and $v_j$ if and only if $v_i\cdot v_j=1$ for $i,
j=1,\ldots,n$. The number of connected components of the intersection
graph of $S$ will be denoted by $c(S)$.

We shall say that an element $v_j\in S$ is~\emph{isolated},~\emph{final}
or~\emph{internal} if the quantity
\[
\sum^n _{\substack{i=1\\ i\neq j}} 
(v_i\cdot v_j) 
\]
is equal to, respectively, $0$, $1$ or $2$. In other words, $v_j$ is
isolated or final if it is, respectively, an isolated vertex or a leaf
of the intersection graph, and it is internal otherwise.

Given elements $e, v\in\bD^n$ with $e\cdot e = -1$, we shall denote by
$\pi_e(v)$ the projection of $v$ in the direction orthogonal to $e$:
\[
\pi_e(v):=v+(v\cdot e) e\in\bD^n.
\]

Two elements $v, w\in\bD^n$ are \emph{linked} if there
exists $e\in\bD^n$ with $e\cdot e=-1$ such that
\[
v\cdot e\neq 0,\quad\text{and}\quad w\cdot e\neq 0.
\]
A set $S\subseteq\bD^n$ is~\emph{irreducible} if, given two elements
$v,w\in S$, there exists a finite sequence
\[
v_0=v, v_1,\ldots, v_k=w \in S
\]
such that $v_i$ and $v_{i+1}$ are linked for  $i=0,\ldots, k-1$. 
A set which is not irreducible is~\emph{reducible}.

The reason to introduce the following definition is technical. It will
become clear later on (see the ``Outline'' at the beginning of
Section~\ref{s:p1>0}).

\begin{defn}\label{d:good}
A subset $S=\{v_1,\ldots,v_n\}\subseteq\bD^n$ is~\emph{good} if it is
irreducible and its elements satisfy~\eqref{e:conds}.
\end{defn}

\begin{defn}\label{d:invariant}
Given a subset $S=\{v_1\ldots,v_n\}\subseteq\bD^n$, let 
\[
I(S):=\sum_{i=1}^n (-v_i\cdot v_i -3) \in \Z
\]
\end{defn}

The following Lemma will be used in
Sections~\ref{s:p1>0},~\ref{s:general}, ~\ref{s:standard}
and~\ref{s:ribbon}.

\begin{lem}\label{l:n=3}
Let $S=\{v_1,v_2,v_3\}\subseteq\bD^3=\langle e_1,e_2,e_3\rangle$ be a
good subset with $I(S)<0$. Then, up to applying to $S$ an element of
$\Aut(\bD^3)$ and possibly replacing $(v_1, v_2,v_3)$ with
$(v_3,v_2,v_1)$, one of the following holds:
\begin{enumerate}
\item
$(v_1,v_2,v_3)=(e_1-e_2,e_2-e_3,-e_2-e_1)$,
\item
$(v_1, v_2,v_3)=(e_1-e_2,e_2-e_3,e_1+e_2+e_3)$,
\item
$(v_1,v_2,v_3)=(e_1+e_2+e_3,-e_1-e_2+e_3,e_1-e_2)$.
\end{enumerate}
Moreover, 
\[
(p_1(S),p_2(S),c(S),I(S))=
\begin{cases}
(1,1,1,-3)\quad\text{in case $(1)$},\\
(0,2,2,-2)\quad\text{in case $(2)$},\\
(0,1,2,-1)\quad\text{in case $(3)$}.
\end{cases}
\]
In particular, $(a_1,a_2,a_3)\in\{(2,2,2),(2,2,3),(3,3,2)\}$.
\end{lem}

\begin{proof}
Up to replacing $(v_1, v_2,v_3)$ with $(v_3,v_2,v_1)$, by
Conditions~\eqref{e:conds} we have three possible cases: (a) $v_1\cdot
v_2=v_2\cdot v_3=1$, (b) $v_1\cdot v_2 = 1$, $v_2\cdot v_3=0$ and (c)
$v_1\cdot v_2=v_2\cdot v_3=0$. Moreover, since $I(S)<0$ we have
$\sum_i a_i \leq 8$. Therefore $a_i\leq 4$ for $i=1,2,3$. Using
the fact that $S$ is irreducible it is easy to see that $a_i<4$ for
$i=1,2,3$. This implies
\[
\{a_1,a_2,a_3\}\in\{\{2,2,2\},\{2,2,3\},\{3,3,2\}\}.
\]
Now observe that if $a_i=3$ then, up to applying an element of
$\Aut(\bD^n)$ we have $v_i=e_1+e_2+e_3$. If $a_j=2$ then $v_j\in\{\pm
e_l\pm e_m\}$, therefore $v_i\cdot v_j$ is an even number, hence
$v_i\cdot v_j=0$. By a similar argument one sees that there cannot be
distinct elements $v_i$ and $v_j$ with $a_i=a_j=3$ and $v_i\cdot
v_j=0$. Using such considerations it is easy to check that, up to
replacing $(v_1, v_2,v_3)$ with $(v_3,v_2,v_1)$,
\begin{enumerate}
\item[(a)]
$(a_1,a_2,a_3)=(2,2,2)$ is the only triple compatible with case (a),
\item[(b)]
$(a_1,a_2,a_3)=(2,2,3)$ is the only triple compatible with case (b),
\item[(c)]
$(a_1,a_2,a_3)=(3,3,2)$ is the only triple compatible with case (c).
\end{enumerate}
The lemma follows by a straightforward case--by--case
analysis.
\end{proof}

The following lemma provides a basic constraint on $p_1(S)$ and
$p_2(S)$ coming from the assumption $I(S)<0$. It will be used in
Sections~\ref{s:p2>0} and~\ref{s:general}.

\begin{lem}\label{l:p1p2ineq}
Let $S\subseteq\bD^n=\langle e_1,\ldots,e_n\rangle$ be a subset of
cardinality $n$ with $I(S)<0$.  Then,
\begin{equation}\label{e:pis}
2p_1(S) + p_2(S) > \sum_{j=4}^n (j-3) p_j(S).
\end{equation}
\end{lem}

\begin{proof}
Let $S=\{v_1,\ldots,v_n\}$ and let $M=(m_{ij})$ be the matrix defined
by $v_i=\sum_j m_{ij} e_j$.  By the definition of $p_i(S)$, the number
of non--zero entries of $M$ is
\[
\sum_{i=1}^n i p_i(S) \leq \sum_{i,j} |m_{ij}| \leq \sum_{i,j} m_{ij}^2 
= - \sum_{i=1}^n v_i\cdot v_i. 
\]
Moreover, the assumption $I(S)<0$ is equivalent to 
\[
- \sum_{i=1}^n v_i\cdot v_i < 3n.
\]
Since it is also evident that 
\[
n = p_1(S) + p_2(S) + \cdots + p_n(S),
\]
the lemma follows.
\end{proof}

Given integers $a_1,\ldots,a_n\geq 2$, we shall use the notation
\[
[a_1,\ldots,a_n]^- := a_1 - \cfrac{1}{a_2 -
       \cfrac{1}{\ddots -
        \cfrac{1}{a_n}
}},
\]
and for any integer $t\geq 0$ we shall write 
\begin{equation}\label{e:power}
(\ldots,2^{[t]},\ldots) := (\ldots,\overbrace{2,\ldots,2}^t,\ldots).
\end{equation}

The following arithmetic lemma will be used in the last section of 
the paper to prove Theorem~\ref{t:main}. 

\begin{lem}\label{l:negsum}
Let $p>q\geq 1$ be coprime integers, and suppose that 
\[
\frac pq = [a_1,\ldots,a_n]^-,\quad
\frac p{p-q} = [b_1,\ldots,b_m]^-,
\]
with $a_1,\ldots,a_n\geq 2$ and $b_1,\ldots,b_m\geq 2$. Then, 
\[
\sum_{i=1}^n (a_i-3) + \sum_{j=1}^m (b_j-3) = -2.
\]
\end{lem}

\begin{proof}
We can write 
\[
\frac pq = 
[m_1,2^{[m_2]},m_3,2^{[m_4]},
\ldots,m_{2s-1},2^{[m_{2s}]}]^-
\]
for some 
\[
m_1,m_3,\ldots,m_{2s-1}\geq 3,\quad
m_2,m_4,\ldots,m_{2s}\geq 0.
\]
Then, by Riemenschneider's point rule~\cite{Ri}
\begin{equation}\label{e:Riemen}
\frac p{p-q} = 
[2^{[m_1-2]},m_2+3,2^{[m_3-3]},
m_4+3,\ldots,2^{[m_{2s-1}-3]},m_{2s}+2]^-.
\end{equation}
Therefore,
\[
\sum_{i=1}^n (a_i-3) = \sum_{i=1}^s (m_{2i-1}-3) - \sum_{i=1}^s m_{2i},
\]
and 
\[
\sum_{j=1}^m (b_j-3) = -1+\sum_{i=1}^s m_{2i} -\sum_{i=1}^s(m_{2i-1}-3) -1.
\]
The lemma follows immediately.
\end{proof}

\section{The case $p_1(S)>0$ and $I(S)<0$}
\label{s:p1>0}

\textbf{Outline.} In this section we introduce the key notion
of~\emph{standard subset}, which is the algebraic object naturally
arising in our approach to Theorem~\ref{t:main} (see the outline of
the proof in Section~\ref{s:intro}). For technical reasons, in order
to understand standard subsets we need to understand first the more
general good subsets introduced in Section~\ref{s:prelim}. In this
section we study the special class of good subsets $S$ satisfying
$p_1(S)>0$ and $I(S)<0$. As explained at the beginning of
Section~\ref{s:p2>0}, this is one of the two important classes of good
subsets $S$ with $I(S)<0$. The main result of this section is
Corollary~\ref{c:p1>0}, which shows that a good subset with $p_1(S)>0$
and $I(S)<0$ is necessarily standard and is obtained from a standard
subset of $\bD^3$ by a finite sequence of operations we
call~\emph{expansions} (see~Definition~\ref{d:contraction}). The results 
of this section will be used in Section~\ref{s:general}. 

\begin{defn}\label{d:standard}
A subset $S_n=\{v_1,\ldots, v_n\}\subseteq\bD^n$
such that
\begin{equation}\label{e:inters0}
v_i\cdot v_j = 
\begin{cases}
-a_i\leq -2\quad\text{if}\quad i=j,\\
1\quad\text{if}\quad |i-j|=1,\\
0\quad\text{if}\quad |i-j|>1.
\end{cases}
\end{equation}
for $i,j=1,\ldots,n$ will be called \emph{standard}.
\end{defn}

The following lemma deals with good subsets $S$ satisfying $p_1(S)>0$.
It will be used in the proofs of Proposition~\ref{p:p1>0},
Corollary~\ref{c:p1>0} and in Section~\ref{s:p2>0}.

\begin{lem}\label{l:p1>0}
Suppose that $n>3$, and let 
\[
S_n=\{v_1,\ldots, v_n\}\subseteq\bD^n=\langle e_1,\ldots,e_n\rangle
\]
be a good subset such that $E^{S_n}_i=\{s\}$ for some
$i,s\in\{1,\ldots,n\}$. Then, 
\begin{enumerate}
\item
$v_s$ is internal;
\item
for some $1\leq j\leq n$ we have $V_s=\{i,j\}$, 
$E_j^{S_n} = \{s-1,s,s+1\}$ and
\[
|v_{s-1}\cdot e_j| =|v_s\cdot e_j|= |v_{s+1}\cdot e_j|=1;
\]
\item
for some $t\in\{s-1,s+1\}$ the set
\[
S_{n-1} := S_n\setminus\{v_s,v_t\}\cup\{\pi_{e_j}(v_t)\}\subseteq
\langle e_1,\ldots,e_{i-1},e_{i+1},\ldots, e_n\rangle\cong\bD^{n-1}
\]
is good, $|E^{S_{n-1}}_j|=1$ and $I(S_{n-1})=I(S_n) + 2 + v_s\cdot
v_s$.
\end{enumerate}
Moreover, if $S_n$ is standard then so is $S_{n-1}$.  
\end{lem}

\begin{proof}
Since $S_n$ is irreducible we have $|V_s|\geq 2$. If $|V_s|>2$, the
set obtained from $S_n$ by replacing $v_s$ with $\pi_{e_i}(v_s)$ would
still satisfy~\eqref{e:conds}, but it would consist of $n$ independent
vectors (see Remark~\ref{r:indep}) contained in the span of the $n-1$
vectors $e_1,\ldots,e_{i-1},e_{i+1},\ldots,e_n$,
giving a contradiction. Therefore $|V_s|=2$, i.e.~$V_s=\{i,j\}$ for
some $j\neq i$. 

If $|v_s\cdot e_j|>1$ then we get a contradiction as before by
replacing $v_s$ with $\pi_{e_i}(v_s)$. Hence, we conclude $|v_s\cdot
e_j|=1$. Since $S_n$ is irreducible and $E^{S_n}_i=\{s\}$, $v_s$ is
not isolated. 

We need to show that $v_s$ is not final. By contradiction, suppose
e.g.~that $v_{s-1}\cdot v_s=0$ and $v_s\cdot v_{s+1}=1$ (the
discussion in the case $v_{s-1}\cdot v_s=1$, $v_s\cdot v_{s+1}=0$ is
similar). Let $l\geq 1$ be the largest natural number such that the
set $\{v_s,\ldots,v_{s+l}\}$ has connected intersection graph. If
\[
a_{s+1}=\cdots=a_{s+l}=2
\]
it is easy to check that $|\cup_{i=0}^l V_{s+i}|=l+2$. Since $S_n$ is
irreducible and $E^{S_n}_i=\{s\}$, this gives a
contradiction. Therefore $a_{s+h}>2$ for some $1\leq h\leq l$. Choose 
$h$ to be as small as possible. Then, it
is easy to verify that for some $k\in\{1,\ldots,n\}$ 
\[
V_{s+h}\cap V_{s+h-1}=\{e_k\}\quad\text{and}
\quad |v_{s+h}\cdot e_k|=1. 
\]
Since $|\cup_{i=0}^{h-1} V_{s+i}|=h+1$, it follows that by eliminating
the vectors 
\[
v_s, v_{s+1},\ldots,v_{s+h-1}
\]
and replacing $v_{s+h}$ with $\pi_{e_k}(v_{s+h})$ one obtains a set of
$n-h$ independent vectors contained in the span of $n-(h+1)$
vectors. This contradiction shows that $v_s$ must be internal,
i.e.~$v_{s-1}\cdot v_s=v_{s+1}\cdot v_s=1$.

Now observe that, since $E^{S_n}_i=\{s\}$, we must have $j\in
V_{s-1}\cap V_{s+1}$.  If $a_{s-1}=a_{s+1}=2$ then $v_{s-1}\cdot
v_{s+1}=0$ implies $V_{s-1}=V_{s+1}$, and it is easy to verify that
either $n=3$ or $S$ is reducible. 

If $a_{s-1},a_{s+1}>2$ then, since clearly $|v_{s-1}\cdot
e_j|=|v_{s+1}\cdot e_j|=1$, one gets a contradiction by eliminating
$v_s$ and replacing $v_{s-1}$ and $v_{s+1}$, respectively, with
$\pi_{e_j}(v_{s-1})$ and $\pi_{e_j}(v_{s+1})$. We conclude that either
(i) $a_{s-1}>2$ and $a_{s+1}=2$ or (ii) $a_{s+1}>2$ and
$a_{s-1}=2$. By symmetry, it suffices to consider the case $a_{s+1}>2$
and $a_{s-1}=2$. Since $|v_{s-1}\cdot e_j|=|v_{s+1}\cdot
e_j|=1$, we have $v_{s-1}\cdot\pi_{e_j}(v_{s+1})=1$. Therefore the
elements of the set
\[
S_{n-1}:=\{v_1,\ldots,v_n\}\setminus\{v_s,v_{s+1}\}
\cup\{\pi_{e_j}(v_{s+1})\}
\]
satisfy~\eqref{e:conds}. Moreover, the formula 
\[
I(S_{n-1})= I(S_n) + 2 + v_s\cdot v_s 
\]
is straightforward to check. Since $E^{S_n}_i=\{s\}$ we have
$E^{S_n}_j=\{s-1,s,s+1\}$, therefore the only vectors linked to $v_s$
are $v_{s-1}$ and $v_{s+1}$. Since $v_{s-1}$ and $\pi_{e_j}(v_{s+1})$
are linked to each other, it follows easily that $S_{n-1}$ is
irreducible. The fact that if $S_n$ is standard then so is $S_{n-1}$
is evident from the definition of $S_{n-1}$.
\end{proof}

The following proposition analyzes the nature of a good subset $S$
with $I(S)<0$ and $p_1(S)>0$. It is essential to prove the main result
of this section, i.e.~Corollary~\ref{c:p1>0}.

\begin{prop}\label{p:p1>0}
Suppose that $n\geq 3$, and let 
\[
S=\{v_1,\ldots, v_n\}\subseteq\bD^n=\langle e_1,\ldots,e_n\rangle
\]
be a good subset such that $I(S)<0$ and $p_1(S)>0$. Then, 
\begin{enumerate}
\item
$S$ is standard;
\item
$|v_i\cdot e_j|\leq 1$ for every $i,j=1,\ldots,n$; 
\item
If $n\geq 4$ there exist $h,t\in\{1,\ldots, n\}$ and $s\in\{1,n\}$ such
that
\[
E_h^S=\{s,t\},\quad a_s=2\quad\text{and}\quad a_t>2.
\]
\end{enumerate}
\end{prop}

\begin{proof}
If $n=3$ the proposition follows from Lemma~\ref{l:n=3}. If $n>3$, set
$S_n:=S$. By Lemma~\ref{l:p1>0} there exists a good subset $S_{n-1}$
with $p_1(S_{n-1})>0$, to which Lemma~\ref{l:p1>0} can be applied
again as long as $n-1>3$. Applying the lemma $n-3$ times we obtain a
sequence $S_n, S_{n-1},\ldots, S_3$ of good subsets with
$p_1(S_n),\ldots,p_1(S_3)>0$. In particular, the fact that $S_3$ is
good and $p_1(S_3)>0$ implies, by Lemma~\ref{l:n=3}, that there is
only one possibility for $S_3$ modulo the action of $\Aut(\bD^3)$,
which is given by Lemma~\ref{l:n=3}(1). This immediately implies that
all the sets $S_i$, $i=3,\ldots,n$, have connected intersection
graph. Therefore $S_n$ is standard, i.e.~(1) holds. Since by
assumption $I(S_n)\leq -1$ and $I(S_3)=-3$, the formula for
$I(S_{n-1})$ in the statement of Lemma~\ref{l:p1>0} implies that every
time we applied the lemma we had $v_s\cdot v_s\geq -4$. Since
$V_s=\{i,j\}$, this implies $v_s\cdot v_s=-2$. Therefore $|v_s\cdot
e_i|=|v_s\cdot e_j|=1$, and by the definition of the sequence $S_n,
S_{n-1},\ldots, S_3$ this immediately implies (3). Finally, it is easy
to check that (2) holds for$S_3$ and $S_4$, and that if $S_{k-1}$ is
obtained from $S_k$ as in Lemma~\ref{l:p1>0} and (2) holds for
$S_{k-1}$ then (2) holds for $S_k$. This proves (2) and concludes the
proof.
\end{proof}

\begin{defn}\label{d:contraction}
Let $S=\{v_1,\ldots,v_n\}\subseteq\bD^n$ be a subset
satisfying~\eqref{e:conds} and such that $|v_i\cdot e_j|\leq 1$ for
every $i,j=1,\ldots,n$. Suppose that there exist $1\leq h,s,t\leq n$
such that
\[
E_h^S=\{s,t\}\quad\text{and}\quad a_t>2.
\]
Then, we say that the subset $S'\subseteq\langle e_1,\ldots,e_{h-1},
e_{h+1},\ldots, e_n\rangle\cong\bD^{n-1}$ defined by
\[
S':=S\setminus\{v_s,v_t\}\cup\{\pi_{e_h}(v_t)\}
\]
is obtained from $S$ by a~\emph{contraction}, and we write $S\searrow
S'$. Moreover, we say that $S$ is obtained from $S'$ by
an~\emph{expansion}, and we write $S'\nearrow S$.
\end{defn}

The following result will be used in the proof of
Corollary~\ref{c:I<0}.

\begin{cor}\label{c:p1>0}
Suppose that $n\geq 3$, and let $S=\{v_1,\ldots, v_n\}\subseteq\bD^n$
be a good subset such that $I(S)<0$ and $p_1(S)>0$. Then, $S$ is standard 
and there is a sequence of expansions
\[
S_3\nearrow S_4\nearrow\cdots\nearrow S_{n-1}\nearrow S_n:=S
\]
such that $S_k$ is standard and $I(S_k)=-3$ for every $k=3,\ldots,n$.
\end{cor}

\begin{proof}
If $n=3$ the corollary follows from Lemma~\ref{l:n=3}. Suppose that
$n\geq 4$, let $S_n:=S$, and let $h, s$ and $t$ be the indexes
appearing in Proposition~\ref{p:p1>0}(3). By Lemma~\ref{l:p1>0} and
Proposition~\ref{p:p1>0}, the set
\[
S_{n-1}:=S_n\setminus\{v_s,v_t\}\cup\{\pi_{e_h}(v_t)\}
\]
is standard and is obtained from $S_n$ by a contraction. Moreover,
$p_1(S_{n-1})>0$, and since $a_s=2$ we have
$I(S_{n-1})=I(S_n)$. Arguing in the same way we get a sequence of
contractions $S_n\searrow S_{n-1}\searrow\cdots\searrow S_3$ with each
$S_k$ standard and $I(S_k)=I(S_3)$ for every $k$. Since by
Lemma~\ref{l:n=3} we have $I(S_3)=-3$, this concludes the proof.
\end{proof}

\section{The case $p_1(S)=0$, $p_2(S)>0$ and $I(S)<0$}
\label{s:p2>0}

\textbf{Outline.}  It follows from Lemma~\ref{l:p1p2ineq} that if
a subset $S\subset\bD^n$ of cardinality $n$ satisfies $I(S)<0$ and
$p_1(S)=0$, then necessarily $p_2(S)>0$. Having dealt with the case
$p_1(S)>0$ in the previous section, in this section we start tackling
the more difficult case of a good subset with $p_1(S)=0$, $p_2(S)>0$
and $I(S)<0$. As in the previous case, one would like to show that
good subsets satisfying this condition are obtained by expansions of
smaller subsets of the same type. But in this case one must first
understand the potential obstruction coming from the fact that during
a sequence of contractions the subset might develop what we
call~\emph{bad components} (see Definition~\ref{d:bad}). The main
result of the section is Proposition~\ref{p:p2>0}, essentially giving
a control on the number of bad components which might appear as a
result of contractions. In the next section we shall use
Proposition~\ref{p:p2>0} to establish some results which hold under
the general assumption $I(S)<0$ and, using these, in
Section~\ref{s:standard} we shall finally be able to show that any
standard subset $S$ with $I(S)<0$ is obtained by expanding standard
subsets of the same type.

\begin{defn}\label{d:bad}
Let $S'=\{v_1,\ldots, v_n\}\subseteq\bD^n$, $n\geq 3$, be a good
subset, and suppose there exists $1<s<n$ such that
$C'=\{v_{s-1},v_s,v_{s+1}\}\subseteq S'$ gives a connected component
of the intersection graph of $S'$, with 
$v_{s-1}\cdot v_{s-1}=v_{s+1}\cdot v_{s+1}= -2$, $v_s\cdot v_s < -2$
and $E^{S'}_j=\{s-1,s,s+1\}$ for some $j$. 
Let $S\subseteq\bD^m$ be a subset of order $m\geq n$ obtained
from $S'$ by a sequence of expansions by final $(-2)$--vectors
attached to $C'$, so that $c(S)=c(S')$ and there is a natural 
$1--1$ correspondence between the sets of connected components of the
intersection graphs of $S$ and $S'$. Then, the connected component
$C\subseteq S$ corresponding to $C'\subseteq S'$ is a~{\bf bad
component} of $S$.  The number of bad components of $S$ will be
denoted by $b(S)$.
\end{defn}

If a good subset $S=\{v_1,\ldots, v_n\}\subseteq\bD^n$ satisfies
$p_2(S)>0$ then for some $i,s,t\in\{1,\ldots,n\}$ we must have
$E^S_i=\{s,t\}$.  There are two possibilities: either $a_s$ and $a_t$
are both greater then $2$, or one of them is equal to $2$. The next
lemma analyzes with the latter possibility (assuming $S$ has no bad
components), while the former possibility is considered in 
Lemma~\ref{l:p2>0II}.

\begin{lem}\label{l:p2>0I}
Suppose that $n>3$, the subset 
\[
S=\{v_1,\ldots, v_n\}\subseteq\bD^n=\langle e_1,\ldots,e_n\rangle
\]
is good, has no bad components and there exist $i,s,t\in\{1,\ldots,n\}$
such that
\[
E^S_i=\{s,t\}\quad\text{and}\quad a_s=2. 
\]
Then, one of the following holds:
\begin{enumerate}
\item 
$v_s\cdot v_t=0$, $v_s$ is not internal, $|V_t|>2$, 
and the set
\[
S':=S\setminus\{v_s,v_t\}\cup\{\pi_{e_i}(v_t)\} 
\subseteq\langle e_1,\dotsc,e_{i-1},e_{i+1}\dotsc, e_n\rangle\cong\bD^{n-1}
\]
is good. Moreover, $I(S')\leq I(S)$ and $S'$ has no bad components.
\item
$v_s\cdot v_t=0$, $v_s$ is internal and $a_t>2$. 
\item
$v_s\cdot v_t=1$, $a_t>2$ and the set $S'$ defined in $(1)$ above is
good. Moreover, $I(S')\leq I(S)$ and $S'$ has no bad components.
\end{enumerate}
\end{lem}

\begin{proof}
Since $a_s=2$, we have $V_s=\{i,j\}$ for some $i,j\in\{1,\ldots,n\}$.

\emph{First case: $v_s\cdot v_t=0$ and $a_t=2$.}

In this case $V_t=\{i,j\}$. Since $S$ is irreducible and $n>3$, there
is a $v_r$ with $r\not\in\{s,t\}$ linked to either $v_s$ or
$v_t$. Since $E_i^S=\{s,t\}$, this implies $j\in V_r$ and $v_r\cdot
v_s=v_r\cdot v_t=1$, therefore $|V_r|\geq 2$. Assuming $|V_r|=2$ it
easily follows that $S$ is reducible. Therefore $|V_r|>2$, which
implies that $S$ has a bad component. Hence this case cannot occur.

\emph{Second case: $v_s\cdot v_t=0$ and $a_t>2$.}

We have $V_t\supseteq\{i,j\}$. Suppose first that $v_s$ is
isolated. If $|V_t|=2$, then no other vector could link $v_s$ nor
$v_t$, and $S$ would be reducible.  If $|V_t|>2$ then the set
\[
S':=S\setminus\{v_s,v_t\}\cup\{\pi_{e_i}(v_t)\} 
\]
satisfies~\eqref{e:conds}. Since every vector $v$ linked to $v_s$ must
satisfy $v\cdot e_j\neq 0$, $S'$ is irreducible. Clearly
$I(S')\leq I(S)$ and it is easy to check that $b(S')=b(S)=0$. Hence
(1) holds.

If $v_s$ is final then $v_s\cdot v_{s'}=1$ for some
$s'\in\{s-1,s+1\}$. This implies $|v_{s'}\cdot e_j|=1$. If $|V_t|>2$
then it follows as above that the set $S'$ is good, $I(S')\leq I(S)$,
and $b(S')=b(S)=0$. Hence (1) holds. If $|V_t|=2$ then $v_{s'}\cdot
v_t=1$, hence $|v_t\cdot e_j|=1$. Since $v_s\cdot v_t=0$, we also have
$|v_t\cdot e_i|=1$. But this is impossible because $a_t>2$.

\emph{Third case: $v_s\cdot v_t=1$ and $a_t=2$.}

In this case $v_s$ is not isolated and $V_t=\{i,k\}$ for some $k\neq
j$. Observe that $v_s$ cannot be a final vector, otherwise
$E^S_j=\{s\}$, which by Lemma~\ref{l:p1>0} implies that $v_s$ is
internal. By symmetry, we may assume without loss of generality that
$t=s+1$. Then, arguing as in the proof of Lemma~\ref{l:p1>0}, one gets
a contradiction using the fact that $E^S_i=\{s,t\}$ by considering the
largest $l,m\geq 1$ such that $\{v_{s-m},\ldots,v_{s+l}\}$ has
connected intersection graph. In fact, if 
\[
a_{s-m}=\cdots = a_{s+l}=2
\]
it is easy to check that $|\cup_{i=s-m}^{s+l} V_i|=l+m+2$. Since $S$
is irreducible and $E_i^S=\{s,s+1\}$, this easily leads to a
contradiction. Therefore, $a_r>2$ for some $s-m\leq r < s$ or
$a_p>2$ for some $s + 1 < p\leq s+l$. Suppose e.g. that only the latter 
happens (the other cases can be dealt with similarly). Choose $p$ as 
small as possible. Then, for some 
$q\in\{1,\ldots,n\}$ 
\[
V_p\cap V_{p-1} = \{e_q\}\quad\text{and}\quad
|v_p\cdot e_q|=1.
\] 
Since $|\cup_{i=s-m}^{p-l} V_i| = p+1-s+m$, one obtains the
contradiction by eliminating the vectors
\[
v_{s-m},\ldots, v_{p-1},
\]
replacing $v_p$ with $\pi_{e_q}(v_p)$ and arguing as in the proof of
Lemma~\ref{l:p1>0}.

\emph{Fourth case: $v_s\cdot v_t=1$ and $a_t>2$.}

By symmetry we may assume $t=s+1$. If $j\not\in V_{s+1}$ then, as in
the case $a_t=2$, $v_s$ is not final, otherwise $E^S_j=\{s\}$, which
implies that $v_s$ is internal by Lemma~\ref{l:p1>0}. Then one obtains
a contradiction as in the previous case by considering the biggest
$l\geq 0$ such that $\{v_{s-l},\ldots,v_s\}$ has connected
intersection graph.

If $V_{s+1}=\{i,j\}$ then $v_s$ cannot be final because otherwise
$v_s$ and $v_{s+1}$ would be linked to no other vector, and therefore
the set $S$ would be reducible. But if $v_s$ is not final then
$E^S_j\supseteq\{s-1,s,s+1\}$ and $k\in V_{s-1}\cap V_{s+1}$ for some
$k\not\in\{i,j\}$, which is impossible if $V_{s+1}=\{i,j\}$.

Therefore we conclude that $V_{s+1}\varsupsetneq\{i,j\}$. Since if
$v_h$ is linked to $v_s$ then $j\in V_h$, it follows that the set
\[
S'=S\setminus\{v_s,v_{s+1}\}\cup\{\pi_{e_i}(v_{s+1})\}
\]
is good, and it is clear that $I(S')\leq I(S)$. Moreover, one can easily
check that $b(S')=b(S)=0$. Hence (3) holds. 
\end{proof}

\begin{lem}\label{l:p2>0II}
Suppose that $n>3$, the subset 
\[
S=\{v_1,\ldots, v_n\}\subseteq\bD^n=\langle e_1,\ldots,e_n\rangle
\]
is good, has no bad components and there exist
$i,s,t\in\{1,\ldots,n\}$ such that
\begin{equation}\label{e:Eis}
E^S_i=\{s,t\}\quad\text{and}\quad 
a_s, a_t > 2.
\end{equation}
Then, up to replacing the pair $(s,t)$ with another pair
satisfying~\eqref{e:Eis}, one of the following holds:
\begin{enumerate}
\item
The set 
\[
S'=S\setminus\{v_s,v_t\}\cup\{\pi_{e_i}(v_t)\}\subseteq
\langle e_1,\dotsc,e_{i-1},e_{i+1},\dotsc,e_n\rangle\cong\bD^{n-1}
\]
is good, $I(S')\leq I(S)-1$ and $b(S')\leq 1$.
\item
There exist $k\neq i$ and $s'\in\{s-1,s+1\}$ such that 
\begin{enumerate}
\item
$E^S_k=\{s,s'\}$,
\item
$v_{s'}\cdot v_s=1$, and 
\item
$a_{s'}=2$.
\end{enumerate}
\end{enumerate}
\end{lem}

\begin{proof}
If the set $S'$ of (1) is good, since $a_s>2$ it follows that
$I(S')\leq I(S)-1$, and it is easy to check that $S'$ can have at most
one bad component, therefore $b(S')\leq 1$. Hence (1) holds.

Now suppose that the set
\[
S'(s,t):=S\setminus\{v_s,v_t\}\cup\{\pi_{e_i}(v_t)\} 
\]
is not good because $\pi_{e_i}(v_t)$ happens to have square $-1$.
Then, $v_t=xe_i\pm e_j$ with $|x|>1$. Since $v_s\cdot v_t\in\{0,1\}$,
we must have $j\in V_s$. Moreover, the vector $\pi_{e_i}(v_s)$ must
have square less than $-1$, because otherwise $v_s=ye_i\pm e_j$ with
$|y|>1$, which implies 
\[
|v_s\cdot v_t|=|-xy\pm 1|\geq 3.
\]
Therefore the set $S'(t,s)$ satisfies~\eqref{e:conds}.  Since there is
no vector linked to $v_t$ but unlinked from $v_s$, it follows that
$S'(t,s)$ is irreducible as well. Therefore, after replacing $(s,t)$
with $(t,s)$, (1) holds.

We may now assume that $S'(s,t)$ and $S'(t,s)$ satisfy~\eqref{e:conds}
but they are not good because they are both reducible. We can write
$S'(s,t)=S'_1\cup S'_2$, where $S'_2$ is a maximal irreducible subset
of $S'(s,t)$ containing $\pi_{e_i}(v_t)$ and $S'_1=S'(s,t)\setminus
S'_2$. Define $S_i\subseteq S\setminus\{v_s\}$, $i=1,2$, to be the
preimage of $S'_i$ under the natural surjective map
$S\setminus\{v_s\}\to S'(s,t)$ sending $v_j$ to $v_j$ for $j\neq s,t$
and $v_t$ to $\pi_{e_i}(v_t)$. The decomposition
$S\setminus\{v_s\}=S_1\cup S_2$ shows that $S\setminus\{v_s\}$ is
reducible. Define
\[
V^j:=\cup_{v_k\in S_j} V_k,\quad j=1, 2.
\]
Observe that $V^1\cap V^2=\emptyset$ and $i\in V^2$, and let
$V_s^j:=V_s\cap V^j$, $j=1, 2$. Since $S$ is irreducible while
$S\setminus\{v_s\}$ is reducible, there exists a vector $v_r\in S_1$
which is linked to $v_s$, therefore $|V_s^1|\geq 1$. If $|V_s^1|>1$
then we could replace $v_s$ with
\[
\tilde v_s := -\sum_{k\in V_s^1} (v_s\cdot e_k) e_k
\]
and $v_t$ with $\pi_{e_i}(v_t)$. For every $u\neq s$ we would have
either $V_u\cap V_s\subseteq V^1_s$ or $V_u\cap V_s\subseteq V^2_s$,
implying that $v_u\cdot\tilde v_s\in\{0,1\}$. The $n$ vectors
resulting from the replacements just described would
satisfy~\eqref{e:conds} and hence be independent, but they would be
contained in the span of $\{e_1,\ldots,e_n\}\setminus\{e_i\}$, giving
a contradiction. Thus, we have $V_s^1=\{k\}$ for some $k$. Since
\[
V^1_s = V_s\cap\bigcup_{v_k\in S_1} V_k = 
\bigcup_{v_k\in S_1} (V_s\cap V_k),
\]
if $v_r\in S_1$ is a vector linked to $v_s$ then $V_r\cap V_s=\{k\}$,
hence $v_r\cdot v_s=1$. This implies that $\{s-1,s+1\}\cap
E^S_k\neq\emptyset$. Moreover, if $v_s$ is final we have
$r\in\{s-1,s+1\}$ and then $E^S_k=\{s-1,s\}$ or $E^S_k=\{s,s+1\}$.  By
symmetry, we may assume that the first case occurs. If $a_{s-1}=2$
then (2) holds. If $a_{s-1}>2$ and $k\not\in V_{s+1}$, we can
eliminate $v_s$, replace $v_{s-1}$ with $\pi_{e_k}(v_{s-1})$ and $v_t$
with $\pi_{e_i}(v_t)$. This gives a contradiction unless
$\pi_{e_k}(v_{s-1})$ happens to have square $-1$. But in that case we
can replace $(s,t)$ with $(s-1,s)$, and by the argument given above
$S'(s-1,s)$ is good, therefore (1) holds.

Now we must consider the case when $v_s$ is not final. We have
$v_{s-1}\cdot v_s=v_s\cdot v_{s+1}=1$, $E^S_k=\{s-1,s,s+1\}$ and
$V_{s-1}\cap V_s=V_s\cap V_{s+1}=\{k\}$. Let us suppose that either
$a_{s-1}>2$ or $a_{s+1}>2$. By symmetry we may assume that
$a_{s-1}>2$. Since $v_t\not\in S_1$, $t\not\in\{s-1,s+1\}$, we have
$v_s\cdot v_t=0$, so we can eliminate $v_s$ and make the
replacements:
\[
v_{s-1}\mapsto \pi_{e_k}(v_{s-1}),\quad 
v_t\mapsto\pi_{e_i}(v_t).
\]
It is easy to see that the resulting set $S''$ of $n-1$ vectors
satisfies~\eqref{e:conds} and can be written as a disjoint union
$S''=S''_1\cup\cdots\cup S''_q$ of maximal irreducible subsets so that
each $S''_l$ is contained in the span of a set of vectors $e_j$ whose
cardinality is equal to $|S''_l|$. We know that for some $l$,
$v_{s+1}\in S''_l$. Moreover, by construction $E^{S''_l}_k=\{s+1\}$.
Since $v_{s+1}$ is not internal in $S''_l$, we get a contradiction
with Lemma~\ref{l:p1>0}(1).

We are left with the case $a_{s-1}=a_{s+1}=2$. In this case it is easy
to deduce that $V_{s-1}=V_{s+1}=\{h,k\}$ for some $h$, and
$E^S_h=\{s-1,s+1\}$. But this means that $S$ contains a bad component,
which is contrary to our assumptions.
\end{proof}

The following is an auxiliary result which will be used in the proof 
of Proposition~\ref{p:p2>0} as well as in Section~\ref{s:general}.

\begin{lem}\label{l:aux}
Let 
\[
S=\{v_1,\ldots,v_n\}\subseteq\bD^n=\langle e_1,\ldots,e_n\rangle 
\]
be a good subset such that $p_1(S)=0$, $p_2(S)>0$ and
$I(S)<0$. Suppose that for each $i,s,t\in\{1,\ldots,n\}$ such that
$E^S_i=\{s,t\}$ we have either $a_s=2$ or $a_t=2$. Then, for at least
one choice of $i,s,t$ such that $a_s=2$ and $E_i^S=\{s,t\}$, either
$v_s$ is not internal or $v_s\cdot v_t=1$.
\end{lem}

\begin{proof}
Suppose by contradiction that for each $i,s,t$ such that $a_s=2$
and $E_i^S=\{s,t\}$ we have $v_s\cdot v_t=0$ and $v_s$ is internal.
Then, if $V_s=\{i,j(i)\}$ it follows immediately that $j(i)\in
V_{s-1}\cap V_s\cap V_{s+1}$, and therefore, since $v_s\cdot v_t=0$
implies $j(i)\in V_t$, we have $|E^S_{j(i)}|\geq 4$. Note that, in 
particular, we must have $n\geq 4$. Consider the
collection $\{j(i)\}$ of all the indices $j(i)$ obtainable in this way. 
Since $p_1(S)=0$ and $|E^S_{j(i)}|\geq 4$ for every $i$, 
by Lemma~\ref{l:p1p2ineq} we have 
\[
p_2(S)>p_4(S)+2p_5(S)+\cdots+(n-3)p_n(S)\geq p_4(S)+p_5(S)+\cdots +p_n(S).
\]
Therefore we must have $j(i)=j(i')$ for some 
$i\neq i'$. But if $E^S_{i'}=\{s',t'\}$ with $a_{s'}=2$ then 
$j(i)=j(i')\in V_{s'-1}\cap V_{s'}\cap V_{s'+1}$ and $j(i)=j(i')\in V_{t'}$, and this easily
implies that $|E^S_{j(i)}|\geq 5$. Moreover, it is easy to check that
if $i, i', i''$ are three distinct indexes with
$|E^S_i|=|E^S_{i'}|=|E^S_{i''}|=2$, then one cannot have
$j(i)=j(i')=j(i'')$. This leads to conflict with
Inequality~\eqref{e:pis}, and we must conclude that there exist
$i,s,t\in\{1,\ldots,n\}$ such that $a_s=2$, $E^S_i=\{s,t\}$ and either
$v_s$ is not internal or $v_s\cdot v_t=1$.
\end{proof}

\begin{prop}\label{p:p2>0}
Suppose that $n>3$ and 
\[
S=\{v_1,\ldots,v_n\}\subseteq\bD^n=\langle e_1,\ldots,e_n\rangle 
\]
is a good subset with no bad components and such that $p_1(S)=0$,
$p_2(S)>0$ and $I(S)<0$. Then, there exist $i,s,t\in\{1,\ldots,n\}$
such that the set
\[
S'=S\setminus\{v_s,v_t\}\cup\{\pi_{e_i}(v_t)\}
\subseteq\langle e_1,\dotsc,e_{i-1},e_{i+1},\dotsc,e_n\rangle
\cong\bD^{n-1}
\]
is good. Moreover, $I(S')\leq I(S)$, $b(S')\leq 1$ and if
$b(S')=1$ then $v_s\cdot v_s < -2$ and  $I(S')\leq I(S)-1$.
\end{prop}

\begin{proof}
Since $p_2(S)>0$, there exist $i,s,t\in\{1,\ldots,n\}$ such that
$E^S_i=\{s,t\}$. If $a_s>2$ and $a_t>2$ the hypotheses of
Lemma~\ref{l:p2>0II} are satisfied. Therefore, since $S$ has no bad
components, the conclusions of Lemma~\ref{l:p2>0II}(1) or
Lemma~\ref{l:p2>0II}(2) hold. 

If the conclusion of Lemma~\ref{l:p2>0II}(1) holds, then the proposition 
follows immediately. If the conclusion of Lemma~\ref{l:p2>0II}(2)
holds then Lemma~\ref{l:p2>0I}(3) applies and (2) holds. Therefore,
from now on we assume that for each $i,s,t\in\{1,\ldots,n\}$ such that
$E^S_i=\{s,t\}$ we have either $a_s=2$ or $a_t=2$. Since $p_1(S)=0$,
$p_2(S)>0$ and $I(S)<0$, by Lemma~\ref{l:aux}, for at least one choice
of $i,s,t$ we have $E^S_i=\{s,t\}$, $a_s=2$ and either $v_s$ is not
internal or $v_s\cdot v_t=1$. Since we are assuming $b(S)=0$, by
Lemma~\ref{l:p2>0I} we see that either the conclusion of
Lemma~\ref{l:p2>0I}(1) or the conclusion of Lemma~\ref{l:p2>0I}(3)
holds.  In both cases the proposition is proved.
\end{proof}

\section{The general case $I(S)<0$}
\label{s:general}

\textbf{Outline.}  In this section we study good subsets $S$ with no
bad components and $I(S)<0$. Our aim is to establish some results
which will be used in the next section to analyze standard subsets
with $I(S)<0$. The main result is Corollary~\ref{c:I<0}, which implies
that a good subset $S$ with no bad components and $I(S)<0$ has
$I(S)\in\{-1,-2,-3\}$ and is obtained by a sequence of expansions from
a subset of $\bD^3$ of the same kind as $S$.  

The following simple lemma is used in the proof of
Proposition~\ref{p:I<0a}.

\begin{lem}\label{l:standnobad}
Let $S\subseteq\bD^n$, $n\geq 3$, be a standard subset with $I(S)<0$. 
Then, $S$ has no bad components. 
\end{lem}

\begin{proof}
Since the intersection graph of a standard subset is connected, if the
only connected component of the intersection graph of $S$ is bad, then
by definition $S$ is obtained via expansions by final $(-2)$--vectors
from a subset $S'\subseteq\bD^3$ consisting of a single bad component
with $I(S')<0$. But Lemma~\ref{l:n=3} forbids the existence of such a
subset.
\end{proof}

The following proposition should be thought of as a generalization of
Proposition~\ref{p:p1>0}(2).

\begin{prop}\label{p:I<0a}
Suppose that $n\geq 3$, and let 
\[
S=\{v_1,\ldots, v_n\}\subseteq\bD^n=\langle e_1,\ldots,e_n\rangle
\]
be a good subset with no bad components such that $I(S)<0$. Then,
$|v_i\cdot e_j|\leq 1$ for every $i,j=1,\ldots,n$.
\end{prop}

\begin{proof}
We argue by induction on $n\geq 3$. If $n=3$ the conclusion is an
immediate consequence of Lemma~\ref{l:n=3}. Therefore, from now on we
shall assume $n>3$. Since $I(S)<0$, by Lemma~\ref{l:p1p2ineq}
Inequality~\eqref{e:pis} holds, therefore either $p_1(S)>0$ or
$p_2(S)>0$. If $p_1(S)>0$ then the conclusion holds by
Proposition~\ref{p:p1>0}, hence we may assume $p_1(S)=0$
and $p_2(S)>0$.

Since $b(S)=0$, by Proposition~\ref{p:p2>0} there is a good subset
$S'\subseteq\bD^{n-1}$ such that
\[
I(S') + b(S') \leq I(S) + b(S) < 0
\]
In particular, it follows that $I(S')<0$. Now we set $S_1:=S$,
$S_2:=S'$, $n_1=n$ and $n_2=n-1$. If $n-1=3$ we stop. If $n-1\geq 4$
and $p_1(S_2)>0$, then by Corollary~\ref{c:p1>0} $S_2$ is standard,
$I(S_2)=-3$ and $S_2$ contracts to a standard subset
$S_3\subseteq\bD^{n-2}$ such that $I(S_3)=-3$. If $n-1\geq 4$ and
$p_1(S_2)=0$ then, since $I(S_2)<0$ we have $p_2(S_2)>0$. If $S_2$ has
a bad component then there is a sequence of contractions from $S_2$ to
a good subset $S'_3$ with a connected component
$\{v_{s-1},v_s,v_{s+1}\}$ such that $v_{s-1}\cdot v_{s-1}=v_{s+1}\cdot
v_{s+1}=-2$, $v_s\cdot v_s < -2$ and $E^{S'_3}_j=\{s-1,s,s+1\}$ for some $j$. 
Then, we set
\[
S_3 = S'_3\setminus\{v_s,v_{s+1}\}\cup\{\pi_{e_j}(v_s)\}.
\]
It is easy to check that $S_3$ is good, has no bad components and
$I(S_3)=I(S_2)+1$. Therefore, in any case we obtain a good subset 
$S_3\subseteq\bD^{n_3}$ with $n_3\geq 2$ and  
\[
I(S_3)+b(S_3)\leq I(S_2)+b(S_2) <0.
\]
Continuing in this way, we obtain a decreasing sequence of good
subsets without bad components 
\[
S_1\subseteq\bD^{n_1},\  S_2\subseteq\bD^{n_2},\cdots, S_k\subseteq\bD^{n_k}
\]
with $n_1>n_2>\cdots>n_k\geq 2$ and 
\[
I(S_{i+1}) + b(S_{i+1}) \leq I(S_i) + b(S_i) < 0,\quad
i=1,\ldots,k-1.
\]
Clearly a good subset $S\subseteq\bD^2$ has $I(S)=-2$. Therefore, by
Lemma~\ref{l:n=3} we have $I(S_k)\geq -3$ and $b(S_k)=0$. Setting
\[
\xi(S_i):=
\begin{cases}
I(S_i)+b(S_i),\quad & i=1,\ldots,k, \\
-3,\quad & i=k+1,
\end{cases} 
\]
since $b(S_1)=0$ we have
\[
\sum_{i=1}^k (\xi(S_i)-\xi(S_{i+1})) = \xi(S_1)-\xi(S_{k+1}) = I(S_1)+3 
\leq -1 + 3 = 2.
\]
Since $\xi(S_i)-\xi(S_{i+1})\geq 0$ for every $i$, we conclude 
\begin{equation}\label{e:bound}
0\leq\xi(S_i)-\xi(S_{i+1}) \leq 2, \quad i=1,\ldots,k.
\end{equation}
In particular, $\xi(S_1)-\xi(S_2)\leq 2$. By a simple calculation 
one easily sees that this is equivalent to 
\[
a_s + |v_t\cdot e_i|^2 \leq b(S_2) + 5. 
\]
\emph{First case: $b(S_2)=0$.} In this case $a_s\geq 2$ implies 
$|v_t\cdot e_i|^2=1$ and therefore $a_s\leq 4$. Since $|V_s|\geq 2$, 
we necessarily have $|v_s\cdot e_j|\leq 1$ for every $j=1,\ldots,n$, 
and this easily implies the statement of the proposition. 

\emph{Second case: $b(S_2)=1$.} Clearly, either $|v_t\cdot e_i|^2=1$
or $|v_t\cdot e_i|^2=4$. If $|v_t\cdot e_i|^2=1$ then $a_s\leq 5$. If
$a_s<5$ the conclusion follows as in the previous case. If $a_s=5$ and
$|v_s\cdot e_k|>1$ for some $k$, then $v_s=ae_i+be_j$ with $a^2+b^2=5$
for some $i,j$. This immediately implies $v_s\cdot v_t=1$, $|v_s\cdot
e_i|=1$ and $V_s\cap V_t=\{i\}$, and one gets a contradiction e.g.~by
replacing $v_s$ with $\pi_{e_i}(v_s)$ and $v_t$ with $\pi_{e_i}(v_t)$.
If $|v_t\cdot e_i|^2=4$ then $a_s=2$. But then
Proposition~\ref{p:p2>0} is incompatible with the assumption
$b(S_2)=1$.
\end{proof}

The following proposition shows that good subsets with no bad
components, possibly disconnected intersection graphs and sufficiently
negative invariant $I(S)$ can be contracted to subsets having the same
properties. This fact will quickly lead to the main result of this
section, i.e.~Corollary~\ref{c:I<0}.

\begin{prop}\label{p:I<0b}
Suppose that $n\geq 4$, and let 
\[
S=\{v_1,\ldots,v_n\}\subseteq\bD^n=\langle e_1,\ldots,e_n\rangle
\]
be a good subset with no bad components such that $I(S)<0$.  Then, for
some $i,s,t$ the set
\[
S'=S\setminus\{v_s,v_t\}\cup\{\pi_{e_i}(v_t)\}
\subseteq\langle e_1,\dotsc,e_{i-1},e_{i+1},\dotsc,e_n\rangle
\cong\bD^{n-1}
\]
is good and has no bad components. Moreover, either
\[
(I(S'),c(S'))=(I(S),c(S))
\]
or
\[
I(S')\leq I(S)-1\quad\text{and}\quad c(S')\leq c(S)+1.
\]
\end{prop}

\begin{proof}
Since $I(S)<0$, by Proposition~\ref{p:I<0a} we have $|v_i\cdot
e_j|\leq 1$ for every $i$ and $j$. Moreover, Inequality~\eqref{e:pis}
holds by Lemma~\ref{l:p1p2ineq}, hence either $p_1(S)>0$ or
$p_2(S)>0$. If $p_1(S)>0$ then the conclusion holds by
Proposition~\ref{p:p1>0}. Therefore from now on we shall assume
$p_1(S)=0$ and $p_2(S)>0$. Under this assumption there exist
$i,s,t\in\{1,\ldots,n\}$ such that $E^S_i=\{s,t\}$. If $a_s>2$ and
$a_t>2$, since $S$ has no bad component the hypotheses of
Lemma~\ref{l:p2>0II} are satisfied and either Lemma~\ref{l:p2>0II}(1)
or Lemma~\ref{l:p2>0II}(2) holds. If Lemma~\ref{l:p2>0II}(2) holds
then so does Lemma~\ref{l:p2>0I}(3).  But this is impossible because
the proof of Lemma~\ref{l:p2>0I} (see the fourth case) shows that if
$V_s=\{i,j\}$ then $V_t\supseteq\{i,j\}$, which is incompatible with
$v_s\cdot v_t=1$ and $|v_t\cdot e_j|\leq 1$ for every $j$. If the
conclusion of Lemma~\ref{l:p2>0II}(1) holds then the set
\[
S'=S\setminus\{v_s,v_t\}\cup\{\pi_{e_i}(v_t)\}
\]
is good and clearly $I(S')\leq I(S)-1$ and $c(S')\leq c(S)+1$.  Now we
argue that $S'$ has no bad components.

First, we claim that if $S'$ has a bad component $C'\subseteq S'$
then, if $\pi\co S\setminus\{v_s\}\to S'$ denotes the natural map,
$v_t\in C:=\pi^{-1}(C')$. Observe that if $v_t\not\in C$ then $v_s$
must be orthogonal to $C$ (i.e. to every element of $C$). Otherwise,
it is easy to check that $v_s$ would have nontrivial intersection with
at least $2$ vectors $e_j$ orthogonal to $S\setminus C$. But then,
adding to $S_1$ the vector obtained from $v_s$ by eliminating all the
vectors $e_j$ which are not orthogonal to $S\setminus C$ would give a
contradiction via rank considerations. We conclude that if $v_t\not\in
C$ then $v_s$ must be orthogonal to $C$. But then if $S'$ has a bad
component also $S$ has one, so we get a contradiction. Therefore
$v_t\in C$.

Next, we observe that by the proof of Proposition~\ref{p:I<0a} we have
$-4\leq v_s\cdot v_s\leq -3$. Using this fact together with $|v_s\cdot
e_j|=1$ for every $j$ it is now a simple exercise to find a
contradiction by analyzing separately the following three cases.
We sketch the argument for each case. 

{\em First case: $v_s$ orthogonal to $C$.}

Since $v_s\cdot v_t=0$, $V_s\cap V_t\supseteq\{i,j\}$ for some
index $j$. Moreover, $|V_s\cap V_t|$ must be even, therefore it is
either $2$ or $4$. But if $|V_s\cap V_t|=4$ one immediately gets a
contradiction from the fact that $v_s$ is orthogonal to $C$.
Therefore, $|V_s\cap V_t|=2$ and $e_j\cdot v'\neq 0$ for some $v'\in
C\setminus\{v_t\}$. Since $E^S_i=\{s,t\}$ and $v_s\cdot v'=0$, for
some $k\neq i,j$ we have $k\in V_s$ and $e_k\cdot v'\neq 0$. It is
now easy to see that $e_k\cdot v''\neq 0$ for some $v''\in
C\setminus\{v_t,v'\}$. Since $v_s\cdot v''=0$,
there is an index $h\neq i,j,k$ such that $k\in V_s$ and $e_k\cdot
v''\neq 0$. Since $|V_s|\leq 4$, continuing in this way clearly leads
to a contradiction. 

{\em Second case: $v_s$ not orthogonal to $C$ but $v_s\cdot v_t=0$.}

As in the previous case, $V_s\cap V_t$ has either $2$ or $4$ elements.
But $|V_s\cap V_t|=4$ easily leads to a contradiction, therefore
$V_s\cap V_t=\{i,j\}$. Let $v_r\in C$ with $v_s\cdot v_r=1$. We have
$V_s\cap V_r\neq\{j\}$ (otherwise a contradiction follows
immediately).  If $V_s=\{i,j,h\}$ then a contradiction follows quickly
by considering the vectors of $C$ which intersect non--trivially
$e_h$. If $V_s=\{i,j,h,k\}$ one gets a contradiction via a rank
counting argument by replacing $v_s$ with $\pi_{e_j}(\pi_{e_i}(v_s))$
and $v_t$ with $\pi_{e_i}(v_t)$.

{\em Third case: $v_s\cdot v_t=1$.}

If $V_s\cap V_t=\{i\}$ then replacing $v_s$ with $\pi_{e_i}(v_s)$ and
$v_t$ with $\pi_{e_i}(v_t)$ one gets a contradiction by the usual rank
counting argument. Therefore $V_s\cap V_t=\{i,j,k\}$. Again, this
gives a contradiction by looking at the vectors of $C$ which intersect
non--trivially $e_j$ and $e_k$. 

The previous arguments show that if $a_s>2$ and $a_t>2$ then the
statement of the proposition holds. Therefore we may now assume that for each
$i,s,t\in\{1,\ldots,n\}$ such that $E^S_i=\{s,t\}$ we have either
$a_s=2$ or $a_t=2$. By Lemma~\ref{l:aux}, for at least one choice of
$i,s,t$ we have $a_s=2$ and either $v_s$ is not internal or $v_s\cdot
v_t=1$. Therefore, since $S$ has no bad component the conclusion of
either Lemma~\ref{l:p2>0I}(1) or Lemma~\ref{l:p2>0I}(3) holds. But, as
we pointed out above, the conclusion of~\ref{l:p2>0I}(3) leads to a
contradiction, therefore~\ref{l:p2>0I}(1) must hold. Thus, the
resulting $S'$ has no bad components and, since $|v_t\cdot e_i|=1$ and
$v_s$ is not internal, we have $I(S')=I(S)$ and $c(S')=c(S)$.
\end{proof}

\begin{cor}\label{c:I<0}
Suppose that $n\geq 3$, and let $S_n=\{v_1,\ldots,v_n\}\subseteq\bD^n$
be a good subset with no bad components and such that $I(S_n) <
0$. Then $I(S_n)\in\{-1,-2,-3\}$, there exists a sequence of
contractions
\begin{equation}\label{e:expans}
S_n\searrow S_{n-1}\searrow\cdots\searrow S_3
\end{equation}
such that, for each $k=3,\dotsc,n-1$ the set $S_k$ is good, has no
bad components and we have either
\[
(I(S_k),c(S_k))=(I(S_{k+1}),c(S_{k+1}))
\]
or
\[
I(S_k)\leq I(S_{k+1})-1\quad\text{and}\quad c(S_k)\leq c(S_{k+1})+1.
\]
Moreover:
\begin{enumerate}
\item
If $p_1(S_n)>0$ then $I(S_n)=-3$, $S_n$ is standard and one can choose
the above sequence so that $I(S_k)=-3$ and $S_k$ is standard for every
$k=3,\ldots, n-1$.
\item
If $I(S_n)+c(S_n)\leq 0$ then $S_3$ is given, up to applying an
automorphism of $\bD^3$, by either (1) or (2) in Lemma~\ref{l:n=3}; if
$I(S_n)+c(S_n)<0$ then the former case occurs.
\end{enumerate}
\end{cor}

\begin{proof}
If $n=3$ the corollary follows immediately from Lemma~\ref{l:n=3}, so
we may assume $n\geq 4$. Since $I(S_n)\leq -c(S_n)<0$, by
Lemma~\ref{l:p1p2ineq} either $p_1(S_n)>0$ or $p_2(S_n)>0$.  If
$p_1(S_n)>0$ then the existence of the required sequence as well as
(1) follow from Corollary~\ref{c:p1>0}. Moreover, in this case (2)
follows from (1) because, by Lemma~\ref{l:n=3}, $S_3$ is given, up to
the action of $\Aut(\bD^3)$, by Lemma~\ref{l:n=3}(1).

If $p_1(S_n)=0$ and $p_2(S_n)>0$ the existence of the
sequence~\eqref{e:expans} follows from several applications of
Proposition~\ref{p:I<0b}.  Since $I(S_k)\leq I(S_{k+1})$ for
$k=3,\ldots,n-1$, $I(S_n)\in\{-1,-2,-3\}$. If $I(S_n)+c(S_n)\leq 0$,
since
\begin{equation}\label{e:ine}
I(S_3)+c(S_3)\leq I(S_4)+c(S_4)\leq \cdots \leq I(S_n)+c(S_n)\leq 0, 
\end{equation}
it follows from Lemma~\ref{l:n=3} that, up to applying an automorphism
of $\bD^3$, $S_3$ must be either of type~\ref{l:n=3}(1)
or~\ref{l:n=3}(2). Inequalities~\eqref{e:ine} imply that if
$I(S_n)+c(S_n)<0$ then $I(S_3)+c(S_3)\leq -1$, hence $S_3$ is given,
up the action of $\Aut(\bD^3)$, by~\ref{l:n=3}(1).
\end{proof}

\section{Standard subsets}
\label{s:standard}

\textbf{Outline.}  In this section we finally look at the subsets of
$\bD^n$ we are mostly interested in, that is the standard subsets with
$I(S)<0$. By Corollary~\ref{c:I<0} such subsets satisfy
$I(S)\in\{-1,-2,-3\}$.  As it turns out, the case $I(S)=-3$ is the
easiest, so we deal with this case first in
Proposition~\ref{p:I=-3}. Theorem~\ref{t:standard} is the main result
and the culmination of all the work done in this section and in the
previous three sections. It is the main algebraic result underlying
the proof of Theorem~\ref{t:main}.  Proposition~\ref{p:I=-3} and
Theorem~\ref{t:standard} will be used in the next section to
characterize the strings $(a_1,\ldots,a_n)$ associated to standard
subsets $S\subset\bD^n$ with $I(S)<0$.

\begin{prop}\label{p:I=-3}
Let $n\geq 3$, and let 
\[
S_n=\{v_1,\ldots,v_n\}\subseteq\bD^n=\langle e_1,\ldots,e_n\rangle 
\]
be a standard subset such that $I(S_n)=-3$. Then, there is a sequence
of contractions
\[
S_n\searrow\cdots\searrow S_3
\]
with $I(S_k)=-3$ and $S_k$ standard for $k=3,\dotsc,n$. Moreover,
\begin{enumerate}
\item
$p_1(S_n)=p_2(S_n)=1$ and $p_3(S_n)=n-2$;
\item
If $E^{S_n}_i=\{s\}$ then $v_s$ is internal and $v_s\cdot v_s=-2$;
\item
If $|E^{S_n}_j|=2$ then $E^{S_n}_j=\{1,n\}$;
\item
either $v_1\cdot v_1=-2$ or $v_n\cdot v_n=-2$. 
\end{enumerate}
\end{prop}

\begin{proof}
We argue by induction on $n$. For $n=3$ the statement of the proposition 
follows immediately from Lemma~\ref{l:n=3}, 
because $I(S_3)=-3$ implies that $S_3$ is given, up to the 
action of $\Aut(\bD^3)$, by~\ref{l:n=3}(1). Let us now assume 
$n>3$. By Corollary~\ref{c:I<0}
there is a sequence of contractions
\[
S_n\searrow\cdots\searrow S_3
\] 
with $I(S_n)\geq \cdots\geq I(S_3)$. Since by Lemma~\ref{l:n=3}
$I(S_3)\geq -3$, the assumption $I(S_n)=-3$ implies $I(S_n)=\cdots
=I(S_3)=-3$. By Corollary~\ref{c:I<0} it follows that each $S_k$ is
standard for $k=3,\ldots,n$. Up to applying an element of
$\Aut(\bD^n)$ we have
\[
S_{n-1} = S_n\setminus\{v_s,v_t\}\cup\{\pi_{e_n}(v_t)\}
\]
for some $s,t$ with $v_s$ final and $v_s\cdot v_s=-2$. Moreover, we
may assume without loss that $v_s=e_1+e_n$. Then, it is easy to check
that $|E^{S_{n-1}}_1|=2$ and therefore, by the induction hypothesis,
$E_1^{S_{n-1}}=\{1,n-1\}$. It follows immediately that $|E_1^{S_n}|=3$
and $E^{S_n}_n=\{1,n\}$, and using this it is very easy to verify the
statement of the proposition for $S_n$.
\end{proof}

The next two lemmas will be used in the proof of Theorem~\ref{t:standard}. 

\begin{lem}\label{l:I=-2}
Let $S_3\nearrow\cdots\nearrow S_n$ be a sequence of expansions such
that, for each $k=3,\ldots,n$, $S_k$ is good, has no bad component and
$(I(S_k),c(S_k))=(-2,2)$. Then,
\begin{enumerate}
\item
$p_1(S_n)=0$, $p_2(S_n)=2$ and $p_3(S_n)=n-2$. 
\item
If $E^{S_n}_i=\{t,t'\}$ then $v_t$ and $v_{t'}$ are not internal and 
exactly one of them has square $-2$. 
\item
If $v_t\in S_n$ is not internal then there exists $i\in V_t$ such that 
$|E^{S_n}_i|=2$.
\end{enumerate}
\end{lem}

\begin{proof}
We argue by induction on $n\geq 3$. For $n=3$ the statement of the lemma 
follows immediately from Lemma~\ref{l:n=3}, because $I(S_3)=-2$
implies that $S_3$ is given, up to the action of $\Aut(\bD^3)$,
by~\ref{l:n=3}(2). Now we assume $n>3$. Up to applying an element of
$\Aut(\bD^n)$ we have
\[
S_{n-1} = S_n\setminus\{v_s,v_t\}\cup\{\pi_{e_n}(v_t)\}
\]
for some $s,t$ with $v_s$ final and $v_s\cdot v_s=-2$. As in the proof
of Proposition~\ref{p:I=-3} we may assume without loss that
$v_s=e_1+e_n$, and it is easy to check that $|E^{S_{n-1}}_1|=2$. 
Using the fact that, by the induction hypothesis, the lemma holds
for $S_{n-1}$ it is now easy to check that $|E_1^{S_n}|=3$ and $v_t$ is
not internal, and from this that the lemma holds for $S_n$.
\end{proof}

\begin{lem}\label{l:extracted}
Let $S_3\subset \bD^3$ be a good subset with $I(S_3)=-3$ and $c(S_3)=1$.
Suppose that $S_3\nearrow\cdots\nearrow S_k$ is a sequence of expansions
such that, for each $h=3,\ldots,k$, $S_h$ is good, has no bad component and
$(I(S_h),c(S_h))=(-3,1)$. Then, it is not possible to expand $S_k$ by 
an isolated $(-3)$--vector. 
\end{lem}

\begin{proof}
We may assume that 
\[
S_k=\{v_1,\ldots,v_k\}\subset\langle e_1,\ldots,e_k\rangle\cong\bD^k.
\]
By contradiction, suppose that $S_{k+1}\subset\bD^{k+1}$ is obtained by expanding $S_{k+1}$ via an isolated $(-3)$--vector $v_{k+1}$. Up to applying an automorphism of $\bD^{k+1}$ we can write $v_{k+1}=e_1+e_2+e_{k+1}$. 
Since $v_{k+1}$ is isolated and $|E^{S_{k+1}}_{k+1}|=2$, we have 
\begin{multline*}
|E_1^{S_{k+1}\setminus\{v_{k+1}\}}| + |E_2^{S_{k+1}\setminus\{v_{k+1}\}}| + |E_{k+1}^{S_{k+1}\setminus\{v_{k+1}\}}| =\\ |E_1^{S_k}| + |E_2^{S_k}| + 1 \equiv \sum_{i\neq k+1} v_i\cdot v_{k+1} \pmod 2 \equiv 0 \pmod 2.
\end{multline*}
This shows that the sum $|E^{S_k}_1| + |E^{S_k}_2|$ must be odd. Therefore by Proposition~\ref{p:I=-3} we may assume $E^{S_k}_1=\{1,k\}$ and either $|E^{S_k}_2|=1$ or $|E^{S_k}_2|=3$. Since $v_{k+1}$ is orthogonal to each $v_i$ for $i=1,\ldots, k$, using Proposition~\ref{p:I=-3} it is easy to check that both cases $|E^{S_k}_2|=1$ and $|E^{S_k}_2|=3$ lead to a contradiction.
\end{proof}

\begin{thm}\label{t:standard}
Let $n\geq 3$, and let 
\[
S_n=\{v_1,\ldots,v_n\}\subseteq\bD^n
\]
be a standard subset such that $I(S_n)<0$. Then,
$I(S_n)\in\{-1,-2,-3\}$ and there is a sequence of contractions
\[
S_n\searrow\cdots\searrow S_3
\]
such that for every $k=3,\ldots,n-1$ the set $S_k$ is standard and
$I(S_k)\leq I(S_{k+1})$.
\end{thm}

\begin{proof}
We argue by induction on $n\geq 3$. For $n=3$ the theorem follows 
immediately from Lemma~\ref{l:n=3}, so we assume
$n>3$ and that the theorem holds true for sets of cardinality
between $3$ and $n-1$. By Corollary~\ref{c:I<0} we have 
$I(S_n)\in\{-1,-2,-3\}$ and there is a sequence of contractions
\[
S_n\searrow\cdots\searrow S_3
\]
such that for every $k=3,\dotsc,n-1$, each $S_k$ is good, it has no
bad components, and we have either
\begin{equation}\label{e:eq1}
(I(S_k),c(S_k))=(I(S_{k+1}),c(S_{k+1}))
\end{equation}
or
\begin{equation}\label{e:eq2}
I(S_k)\leq I(S_{k+1})-1\quad\text{and}\quad c(S_k)\leq c(S_{k+1})+1.
\end{equation}
If $I(S_n)=-3$ the theorem follows from Proposition~\ref{p:I=-3},
therefore we may assume $I(S_n)\in\{-2,-1\}$.

Suppose first $I(S_n)=-2$. By Corollary~\ref{c:I<0}(2) we have
$(I(S_3),c(S_3))=(-3,1)$. Then, Equations~\eqref{e:eq1}
and~\eqref{e:eq2} force $c(S_k)=1$ for every $k=3,\ldots, n-1$,
therefore the theorem follows in this case.

Now assume $I(S_n)=-1$. By Equations~\eqref{e:eq1} and~\eqref{e:eq2}
we have $c(S_{n-1})\leq 2$. If $c(S_{n-1})=1$ we can apply the
induction hypothesis and immediately obtain the theorem, therefore
we may assume $(I(S_{n-1}),c(S_{n-1}))=(-2,2)$. By
Corollary~\ref{c:I<0}(2) $(I(S_3),c(S_3))$ is equal to either $(-2,2)$
or $(-3,1)$. If $(I(S_3),c(S_3))=(-3,1)$, it is easy to check using~\eqref{e:eq1} and~\eqref{e:eq2} that for some $3\leq k < n-1$ we have $S_{k+1}=\{v_1,\ldots,v_{k+1}\}$ and
\[
S_k = S_{k+1}\setminus\{v_{k+1},v_t\}\cup\{\pi_{e_{k+1}}(v_t)\}, 
\]
where $I(S_{k+1})=-2$ and $I(S_k)=-3$. But again by~\eqref{e:eq1} and~\eqref{e:eq2}
we must have 
\[
I(S_k)=I(S_{k-1})=\cdots = I(S_3)=-3
\]
and therefore $c(S_k)=c(S_{k-1})=\cdots =c(S_3)=1$. This implies that 
$v_{k+1}$ is isolated and $v_{k+1}\cdot v_{k+1}=-3$, but it contradicts Lemma~\ref{l:extracted}. Therefore from now on we assume $(I(S_3),c(S_3))=(-2,2)$.

The contraction $S_n\searrow S_{n-1}$ involves eliminating an internal
vector of square $-3$, while the sequence of contractions
\begin{equation}\label{e:expans2}
S_{n-1}\searrow\cdots\searrow S_3
\end{equation}
satisfies the assumptions of Lemma~\ref{l:I=-2}. Let us write
\[
S_{n-1}=S_n\setminus\{v_s,v_t\}\cup\{\pi_{e_i}(v_t)\}
\subseteq\langle e_1,\ldots,e_{i-1},e_{i+1},\ldots,e_n\rangle
\cong\bD^{n-1}
\]
for some $i$, $1<s<n$ and $t\neq s$ with $a_s=3$. Up to applying an
automorphism of $\bD^n$ we may assume $i=n$ and $v_s=e_1+e_2+e_n$.
Moreover, we can write $S_{n-1}$ as a union $S_{n-1}=S^1_{n-1}\cup
S^2_{n-1}$ of subsets with connected intersection graphs, where
\[
S^1_{n-1}=\{v'_1,\ldots,v'_{s-1}\},\quad
S^2_{n-1}=\{v'_{s+1},\ldots,v'_n\}.
\]
In view of Proposition~\ref{p:I<0a} it is easy to check that, since
$v_{s-1}\cdot v_s=v_{s+1}\cdot v_s=1$ and $|E^{S_n}_n|=2$, we have 
\begin{multline*}
2 = v_{s-1}\cdot v_s + v_{s+1}\cdot v_s = \sum_{i\neq s} v_i\cdot v_s 
\equiv \\ |E^{S_n\setminus\{v_s\}}_1| + |E^{S_n\setminus\{v_s\}}_2| 
+ |E^{S_n\setminus\{v_s\}}_n| \pmod 2 
\equiv |E^{S_{n-1}}_1| + |E^{S_{n-1}}_2| +1 \pmod 2,
\end{multline*}
and therefore the sum $|E^{S_{n-1}}_1|+|E^{S_{n-1}}_2|$ must be
odd. On the other hand, by Lemma~\ref{l:I=-2} $|E^{S_{n-1}}_i|$ is
equal to either $2$ or $3$ for every $i$. Therefore, we may assume
$|E^{S_{n-1}}_1|=2$ and $|E^{S_{n-1}}_2|=3$. Moreover, by
Lemma~\ref{l:I=-2} we may also assume $|E^{S_{n-1}}_3|=2$ and
$E^{S_{n-1}}_1\cup E^{S_{n-1}}_3 = \{1,s-1,s+1,n\}$. Therefore, if
$e_1\in V_{s-1}\cap V_{s+1}$ then $e_3\in V_1\cap V_n$. In this case,
by Lemma~\ref{l:I=-2} either $v_1\cdot v_1=-2$ and $v_n\cdot v_n<-2$
or $v_n\cdot v_n=-2$ and $v_1\cdot v_1<-2$. By symmetry we may assume
the latter occurs and define
\[
S'_{n-1}:=S_n\setminus\{v_1,v_n\}\cup\{\pi_{e_3}(v_1)\}.
\] 
Clearly $(I(S'_{n-1}),c(S'_{n-1}))=(-1,1)$ and $S'_{n-1}$ is obtained
from $S_n$ by a contraction, hence applying the induction hypothesis
to $S'_{n-1}$ we get the statement of the theorem. Thus, by symmetry and
Lemma~\ref{l:I=-2} we may assume $e_1\in V_{s-1}\cap V_n$. If $e_n\in
V_{s+1}$ then, since $v_s\cdot v_{s+1}=1$, $e_2\not\in V_{s+1}\cup
V_{s-1}$. But this conflicts with $|E^{S_{n-1}}_1|=2$ and
$|E^{S_{n-1}}_2|=3$, therefore $e_n\not\in V_{s+1}$, $e_2\in
V_{s+1}$ and $e_3\in V_1\cap V_{s+1}$. If $v_1\cdot v_1 = -2$ then 
by Lemma~\ref{l:I=-2} $v_{s+1}\cdot v_{s+1} < -2$, we can define 
\[
S'_{n-1}:=S_n\setminus\{v_1,v_{s+1}\}\cup\{\pi_{e_3}(v_{s+1})\}
\] 
and argue as before. If $v_1\cdot v_1 < -2$ then 
$v_{s+1}\cdot v_{s+1} = -2$ and therefore $V_{s+1}=\{e_2,e_3\}$. 
Since $v_1\cdot v_{s+1}=0$ and $e_3\in V_1$, this implies $e_2\in V_1$. 
Now either $s=2$ and $v_1=v_{s-1}$ or $s>2$ and $v_1\cdot v_{s-1}=0$. 
In the former case $e_1\in V_1=V_{s-1}$, and since $v_1\cdot v_s=1$, we 
must also have $e_n\in V_1$. In the latter case we still have $e_n\in V_1$ because $e_1\not\in V_1$. Therefore in either case we can define 
\[
S'_{n-1}:=S_n\setminus\{v_1,v_s\}\cup\{\pi_{e_n}(v_1)\}.
\]
and conclude as before. 
\end{proof}

\section{Strings associated to standard subsets}
\label{s:strings}

\textbf{Outline.}  In this section we use Proposition~\ref{p:I=-3} and
Theorem~\ref{t:standard} to identify the strings $(a_1,\ldots,a_n)$
corresponding to standard subsets $S\subseteq\bD^n$ with
$I(S)\in\{-1,-2,-3\}$. These results will be used in
Section~\ref{s:final} to prove Theorem~\ref{t:main}.

\sh{The case $I=-3$}

Recall Notation~\eqref{e:power}.

\begin{lem}\label{l:I=-3(2)}
Let $n\geq 3$ and let $S_n=\{v_1,\ldots,v_n\}\subseteq\bD^n$ be a
standard subset such that $I(S_n)=-3$. Suppose $v_i\cdot v_i=-a_i$ for
$i=1,\ldots,n$. Then, the string $(a_1,\ldots,a_n)$ is obtained from
$(2,2,2)$ via a finite sequence of operations of the following types:
\begin{enumerate}
\item
$(s_1,s_2,\ldots,s_{k-1},s_k)\mapsto (s_1+1,s_2,\ldots,s_{k-1},s_k,2)$,
\item
$(s_1,s_2,\ldots,s_{k-1},s_k)\mapsto
(2,s_1,s_2,\ldots,s_{k-1},s_k+1)$.
\end{enumerate}
It follows that either $(a_1,\ldots,a_n)$ or $(a_n,\ldots,a_1)$ is of
the form
\[
\begin{cases}
{\scriptstyle (c_k+1,2^{[c_{k-1}-1]},c_{k-2}+2,\ldots,c_3+2,2^{[c_2-1]},c_1+2,
2^{[c_1+1]},c_2+2,\ldots,c_{k-1}+2,2^{[c_k-1]})}\quad\text{or}\\ 
(c_1+1,2^{[c_1+1]}),
\end{cases}
\] 
for some integers $c_1,\ldots,c_k\geq 1$ and $k\geq 3$.
\end{lem}

\begin{proof}
By Proposition~\ref{p:I=-3} there is a sequence of expansions
\[
S_3\nearrow\cdots\nearrow S_n
\]
with $S_3$ given, up to applying an element of $\Aut(\bD^3)$, by
Lemma~\ref{l:n=3}(1) and each expansion is obtained by adding a final
vector of square $-2$ while simultaneously decreasing by $1$ the
square of the opposite final vector. This immediately implies the
first part statement. The second part of the statement follows from a
straightforward calculation.
\end{proof}

\sh{The case $I=-2$}

\begin{lem}\label{l:I=-2(2)}
Let $n\geq 4$, and let $S_n=\{v_1,\ldots,v_n\}\subseteq\bD^n$ be a
standard subset such that $I(S_n)=-2$. Suppose $v_i\cdot v_i=-a_i$ for
$i=1,\ldots,n$. Then, either $(a_1,\ldots,a_n)$ or $(a_n,\ldots,a_1)$ is 
of one of the following types:
\begin{enumerate}
\item
$(2^{[t]},3,2+s,2+t,3,2^{[s]})$, $s,t\geq 0$,
\item
$(2^{[t]},3+s,2,2+t,3,2^{[s]})$, $s,t\geq 0$.
\end{enumerate}
\end{lem}

\begin{proof}
By Theorem~\ref{t:standard} and Lemma~\ref{l:n=3} there is a
sequence of contractions
\[
S_n\searrow\cdots\searrow S_k\searrow S_{k-1}\searrow\cdots\searrow S_3
\]
of standard subsets with $I(S_3)=-3$ and therefore, for some $n\geq
k>3$, $I(S_n)=\cdots=I(S_k)=-2$ and $I(S_{k-1})=\cdots = I(S_3)=-3$.
Moreover, we may assume $S_k=\{v_1,\ldots,v_k\}$ and
\[
S_{k-1}=S_k\setminus\{v_k,v_t\}\cup\{\pi_{e_k}(v_t)\}
\subseteq\langle e_1,\ldots,e_{k-1}\rangle\cong
\bD^{k-1}
\]
for some $1\leq t\leq k-1$ and $v_k\cdot v_k=-3$, with $v_k$ final. Up
to applying an element of $\Aut(\bD^k)$ we may also assume that
$v_k=e_1+e_2+e_k$.  Moreover, by Proposition~\ref{p:I=-3} we have
$p_1(S_{k-1})=1$.  Therefore, if $|E^{S_{k-1}}_1|>1$ and
$|E^{S_{k-1}}_2|>1$ then we would have $p_1(S_k)>0$, which would imply
$I(S_k)=-3$ by Corollary~\ref{c:I<0}(1). Since $v_k$ is final, it is
easy to see, as in the proof of Theorem~\ref{t:standard}, that the
number $|E^{S_{k-1}}_1|+|E^{S_{k-1}}_2|$ must be even. Thus, in view
of Proposition~\ref{p:I=-3} we may assume $|E^{S_{k-1}}_1|=1$ and
$|E^{S_{k-1}}_2|=3$. By (2) of the same proposition there is a vector
$v_h\in S_k$ such that $1<h<k-1$ and $e_1\cdot
v_h\neq 0$. If $t=h$ then $v_h\cdot e_k\neq 0$, therefore
$E_1^{S_k}=E_k^{S_k}=\{h,k\}$. Since $|E_2^{S_k}|=4$, there exists
$s\not\in\{h,k-1,k\}$ such that $e_2\in V_s$. But $v_s\cdot v_k=0$
implies $s\in E_1^{S_k}\cup E_k^{S_k}$, which is impossible. Therefore
we have $t\neq h$, $v_h\cdot e_k=0$ and $e_2\in V_h$. Then
$v_h\in\{\pm(e_1-e_2)\}$, and since $|E_1^{S_{k-1}}|=1$ this implies
$e_2\in V_{h-1}\cap V_{h+1}$ as well. If $h+1<k-1$ then, since
$E_1^{S_k}=\{h,k\}$, $v_{h-1}\cdot v_k=v_{h+1}\cdot v_k=0$ implies
$e_k\in E_{h-1}^{S_k}\cap E_{h+1}^{S_k}$, which is impossible because
$|E_k^{S_k}|=2$. Therefore we must conclude $h+1=k-1$ and $e_k\in
V_{h-1}$.

Combining this analysis with the proof of Proposition~\ref{p:I=-3}
shows that if $S_{k-1}=\{v'_1,\ldots,v'_{k-1}\}$ then, 
up to the action of $\Aut(\bD^{k-1})$ we have 
\[
v'_{k-1} = -e_2 - e_3 -\cdots - e_{k-1}\quad\text{and}\quad
v'_{k-i} = e_i-e_{i-1},\ i=2,\ldots,k-1.
\]
Moreover, $v_i=v'_i$ for $1\leq i\leq k-1$, $i\neq k-3$, and
$v_{k-3}=v'_{k-3}+e_k$. Now we see that $S_k$ can be contracted to a
standard subset $S'_{k-1}$ by dropping $v_1$ and replacing $v_{k-1}$
with $v^{(1)}_{k-1}:=\pi_{e_{k-1}}(v_{k-1})$. Similarly, for
$i=2,\ldots,k-4$ we can define $S'_{k-i}$ by dropping $v_i$ from
$S'_{k-i+1}$ and replacing $v^{(i-1)}_{k-i+1}$ with
$v^{(i)}_{k-i}:=\pi_{e_{k-i}}(v^{(i-1)}_{k-i+1})$. Continuing in this
way we can construct a new sequence of contractions
\begin{equation}\label{e:new-contractions}
S_k\searrow S'_{k-1}\searrow\cdots\searrow S'_4
\end{equation}
of standard subsets with $I(S'_i)=-2$ for $i=4,\ldots,k-1$, and such
that, up to an automorphism of $\bD^4$,
\[
S'_4=\{e_3-e_2+e_4,e_2-e_1,-e_2-e_3,e_1+e_2+e_4\}.
\]
Analysing Sequence~\eqref{e:new-contractions} it is easy to see
that, up to reversing the $k$--tuple $(v_1,\ldots,v_k)$, if
$|E_i^{S_k}|=2$ then
\[
E_i^{S_k}\in\{\{1,k-1\},\{k-2,k\},\{k-3,k\}\}. 
\]
Since the subset $S_n$ is obtained from $S_k$ by a sequence of
expansions by final $(-2)$--vectors, the lemma follows easily.
\end{proof}

\sh{The case $I=-1$}

\begin{lem}\label{l:I=-1}
Let $n\geq 4$ and let $S_n=\{v_1,\ldots,v_n\}\subseteq\bD^n$ be a
standard subset such that $I(S_n)=-1$. Suppose $v_i\cdot v_i=-a_i$ for
$i=1,\ldots,n$. Then, either $(a_1,\ldots,a_n)$ or $(a_n,\ldots,a_1)$
is of one of the following types:
\begin{enumerate}
\item
$(t+2,s+2,3,2^{[t]},4,2^{[s]})$, $s,t\geq 0$, 
\item
$(t+2,2,3+s,2^{[t]},4,2^{[s]})$, $s,t\geq 0$,
\item
$(3+t,2,3+s,3,2^{[t]},3,2^{[s]})$, $s,t\geq 0$.
\end{enumerate}
\end{lem}

\begin{proof}
By Theorem~\ref{t:standard} and Lemma~\ref{l:n=3} there is a sequence
of contractions
\[
S_n\searrow\cdots\searrow S_3
\]
of standard subsets with $I(S_3)=-3$. Thus, either for some $3< k\leq
n$ we have
\begin{equation}\label{e:seq2}
I(S_n)=\cdots=I(S_k)=-1,\quad
I(S_{k-1})=\cdots = I(S_3)=-3,
\end{equation}
or for some $3\leq k<h\leq n$ we have 
\begin{equation}\label{e:seq1}
\begin{split}
I(S_n)=\cdots=I(S_h)=-1,\quad
I(S_{h-1})=\cdots=I(S_k)=-2,\\
I(S_{k-1})=\cdots = I(S_3)=-3.
\end{split}
\end{equation}

{\em First case: \eqref{e:seq2} holds}. 

The expansion $S_{k-1}\nearrow S_k$ is obtained by adding a final
$(-4)$--vector which can be assumed of the form $v_k=e_1+e_2+e_3+e_k$,
and $p_1(S_k)=0$ otherwise by Corollary~\ref{c:I<0}(1)
$I(S_k)=-3$. Moreover, by Proposition~\ref{p:I=-3} we have
$p_1(S_{k-1})=1$ and, by the parity argument used in the proofs of
Theorem~\ref{t:standard} and Lemma~\ref{l:I=-2(2)} we have
$|E^{S_{k-1}}_1|=1$, $|E^{S_{k-1}}_2|=2$ and $|E^{S_{k-1}}_3|=3$ (up
to renaming the vectors $e_1$, $e_2$ and $e_3$). Also, by
Proposition~\ref{p:I=-3} we have $e_1\in V_l$ for some $1<l<k-1$ and
$v_l$ is of the form $v_l=\pm e_1 \pm e_i$ with
$E^{S_{k-1}}_i=\{l-1,l,l+1\}$.  Since $v_l\cdot v_k=0$ this
immediately implies $i=3$. Again by the proposition we have $e_2\in
V_1\cap V_{k-1}\cap V_k$.  Since $v_k\cdot v_1=v_k\cdot v_{l-1}=0$, it
is easy to check that if $l-1\neq 1$ then $e_k\in V_{l-1}\cap V_1$,
which is impossible because $|E_k^{S_k}|=2$. Therefore $l=2$, $e_3\in
V_1$ and $e_k\in V_3$. By the proposition this implies that $S_k$ has
associated string (up to a reflection) of the form
\[
(t+2,2,3,2^{[t]},4),
\]
$E^{S_k}_1=\{2,k\}$, $E^{S_k}_k=\{3,k\}$ and $|E_i^{S_k}|>2$ for
$i\neq 1,k$. Since $S_n$ is obtained from $S_k$ by a sequence of
expansions obtained by adding final vectors of square $-2$, this
implies that $S_n$ has associated fraction (up to a reflection) as in
(1) or (2) from the statement of the lemma.

{\em Second case: \eqref{e:seq1} holds.}

Arguing as in the proof of 
Lemma~\ref{l:I=-2(2)} we may assume that $k=4$, and 
\[
S_4=\{v_1=e_1-e_2+e_4,v_2=e_2+e_3,v_3=-e_2-e_1,v_4=e_4+e_2-e_3\},
\]
with $S_4\nearrow\cdots\nearrow S_{h-1}$ consisting of expansions
obtained by adding final $(-2)$--vectors and the
expansion $S_{h-1}\nearrow S_h$ obtained by adding a final
$(-3)$--vector which we can assume to be $v_h=e_1+e_2+e_h$. 
Since the number $|E^{S_{h-1}}_1|+|E^{S_{h-1}}_2|$ must be even and it
can be easily checked that $p_2(S_{h-1})=3$, $p_4(S_{h-1})=1$ and
$p_3(S_{h-1})=h-5$, a case--by--case analysis shows that
$|E^{S_{h-1}}_1|=|E^{S_{h-1}}_2|=2$.  

This implies, assuming $S_h=\{v_1,\ldots,v_h\}$, that $i\in
E_1^{S_h}\cap E_2^{S_h}$ for some $i$ with $i<h$. The same analysis of
the sequence $S_4\nearrow\cdots\nearrow S_{h-1}$ used at the end of
the proof of Lemma~\ref{l:I=-2(2)} shows that $|E_i^{S_{h-1}}|=2$
implies $E_i^{S_{h-1}}\cap\{1,h-1\}\neq\emptyset$. Since $v_{h-1}\cdot
v_h=1$, up to renaming $e_1$ and $e_2$ we may assume $e_1, e_2\in V_1$
and $e_2\in V_{h-1}$. It is easy to check that this implies that all
the $(-2)$--vectors added in the sequence $S_4\nearrow\cdots\nearrow
S_{h-1}$ are added from the same side. If they are added e.g.~from the
right--hand side the string associated to $S_{h-1}$ has the form
\[
(3+t,2,2,3,2^{[t]}).
\]
Moreover, the same analysis as above shows that $E^{S_h}_h=\{3,h\}$
and $|E^{S_h}_i|=3$ for every $i\in V_1$. This implies that the
sequence $S_h\nearrow\cdots\nearrow S_n$ consists of expansions
obtained by adding $(-2)$--vectors from the right--hand side only. It
follows that the string associated to $S_n$ is of the form (3) from
the statement of the lemma. If the $(-2)$--vectors added in the sequence
$S_4\nearrow\cdots\nearrow S_{h-1}$ are added from the left--hand
side the resulting string is obtained from (3) by a reflection.
\end{proof}

\section{Existence of ribbon surfaces}
\label{s:ribbon}

\textbf{Outline.} In this section we prove the existence of bounding
ribbon surfaces for all the 2--bridge links which will occur in the
proof of Theorem~\ref{t:main} in Section~\ref{s:final}.

We shall use the following elementary fact about continued fractions
(see~e.g.~\cite[Appendix]{OW} for a proof). Let $p>q\geq 1$ be coprime
integers, and suppose that
\[
\frac pq=[a_1,a_2,\ldots,a_n]^-,\quad a_1,\ldots,a_h\geq 2.
\]
Then, 
\begin{equation}\label{e:inversion}
\frac p{q'} = [a_n,a_{n-1},\ldots,a_1]^-,
\end{equation}
where $p>{q'}\geq 1$ and $q{q'}\equiv 1\pmod p$. 

Let $a_1,\ldots,a_{2n}$ be positive integers. The following 
identity holds (see~\cite[Proposition~2.3]{PP}):
\begin{equation}\label{e:cfrac}
{\textstyle
[a_1,\ldots,a_{2n}]^+ = 
\begin{cases}
{\scriptstyle
[a_1+1,2^{[a_2-1]},a_3+2,2^{[a_4-1]},\ldots,a_{2n-1}+2,2^{[a_{2n}-1]}]^-},
\quad n\geq 2,\\
[a_1+1,2^{[a_2-1]}]^-,\quad n=1.
\end{cases}
}
\end{equation}

\begin{lem}\label{l:familyI}
Let $p>q\geq 1$ be coprime integers, and suppose that\\ 
$\frac pq = [a_1,\ldots,a_n]^-$, where either $(a_1,\ldots,a_n)$ or
$(a_n,\ldots,a_1)$ is of the form
\[
\begin{cases}
{\scriptstyle (c_k+1,2^{[c_{k-1}-1]},c_{k-2}+2,\ldots,c_3+2,2^{[c_2-1]},c_1+2,
2^{[c_1+1]},c_2+2,\ldots,c_{k-1}+2,2^{[c_k-1]})}\quad\text{or}\\ 
(c_1+1,2^{[c_1+1]}),
\end{cases}
\] 
for some integers $c_1,\ldots,c_k\geq 1$ and $k\geq 3$.  Then, if
$p$ is odd $K(p,q)$ bounds an immersed ribbon disk; if $p$ is even the
2--component link $K(p,q)$ bounds the image under a ribbon immersion
of the disjoint union of a disk and a M\"obius band.
\end{lem}

\begin{proof}
Let $k,c_1,\ldots,c_k\geq 1$ be integers. Then, a straightforward
application of Equation~\eqref{e:cfrac} gives 
\begin{multline*}
{\scriptstyle
[c_k+1,2^{[c_{k-1}-1]},c_{k-2}+2,\ldots,c_3+2,2^{[c_2-1]},c_1+2,
2^{[c_1+1]},c_2+2,\ldots,c_{k-1}+2,2^{[c_k-1]})]^-} = \\
[c_k,c_{k-1},c_{k-2},\ldots,c_1,c_1+2,c_2,c_3,
\ldots,c_{k-1},c_k]^+
\end{multline*}
and
\[
[c_1+1,2^{[c_1+1]}]^- = [c_1,c_1+2]^+.
\]
Recalling that if $0<q'<p$ and $qq'\equiv 1\pmod p$ the link $K(p,q')$
is isotopic to $K(p,q)$, this shows that the $K(p,q)$ is isotopic to
$L(c_k,\ldots,c_1,c_1+2,c_2,\ldots,c_k)$ (see Figure~\ref{f:fig1},
case $n$ even). Applying the ribbon move described in the top picture
of Figure~\ref{f:fig2} reduces $K(p,q)$ to a 2--component unlink, as
shown in the remaining pictures of Figure~\ref{f:fig2}. By standard
facts on ribbon moves, this proves the lemma.
\begin{figure}[ht]
\begin{center}
\psfrag{-c1}{${\scriptstyle -c_1}$} 
\psfrag{c2}{${\scriptstyle c_2}$}
\psfrag{-c3}{${\scriptstyle -c_3}$} 
\psfrag{c1+2}{${\scriptstyle c_1+2}$} 
\psfrag{-c2}{${\scriptstyle -c_2}$} 
\psfrag{ck}{${\scriptstyle c_k}$}
\psfrag{-ck}{${\scriptstyle -c_k}$}
\psfrag{iso}{(isotopy)}
\includegraphics[width=7cm]{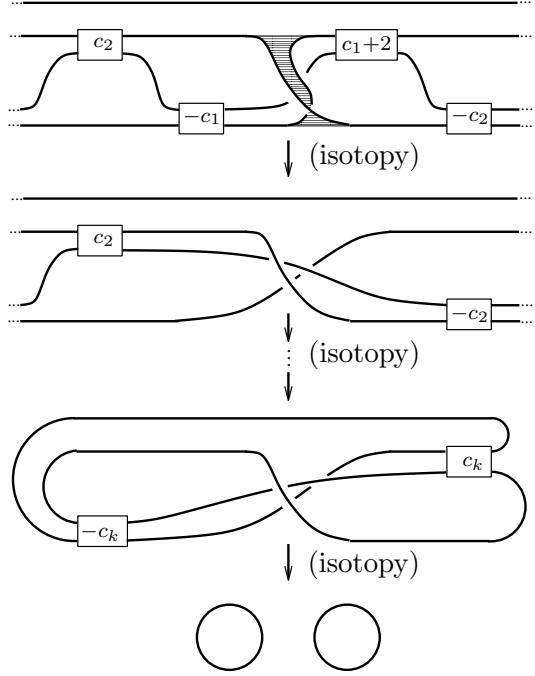}
\end{center}
\caption{The link $L(c_k,\ldots,c_1,c_1+2,\ldots,c_k)$ bounds a
ribbon surface}
\label{f:fig2}
\end{figure}
\end{proof}

\begin{figure}[ht]
\begin{center}
\psfrag{a}{$\scriptstyle a$}
\psfrag{b}{$\scriptstyle b$}
\psfrag{-a}{$\scriptstyle -a$}
\psfrag{-b}{$\scriptstyle -b$}
\psfrag{ribbon}{(two ribbon moves)}
\psfrag{iso}{(isotopy)}
\includegraphics[width=8cm]{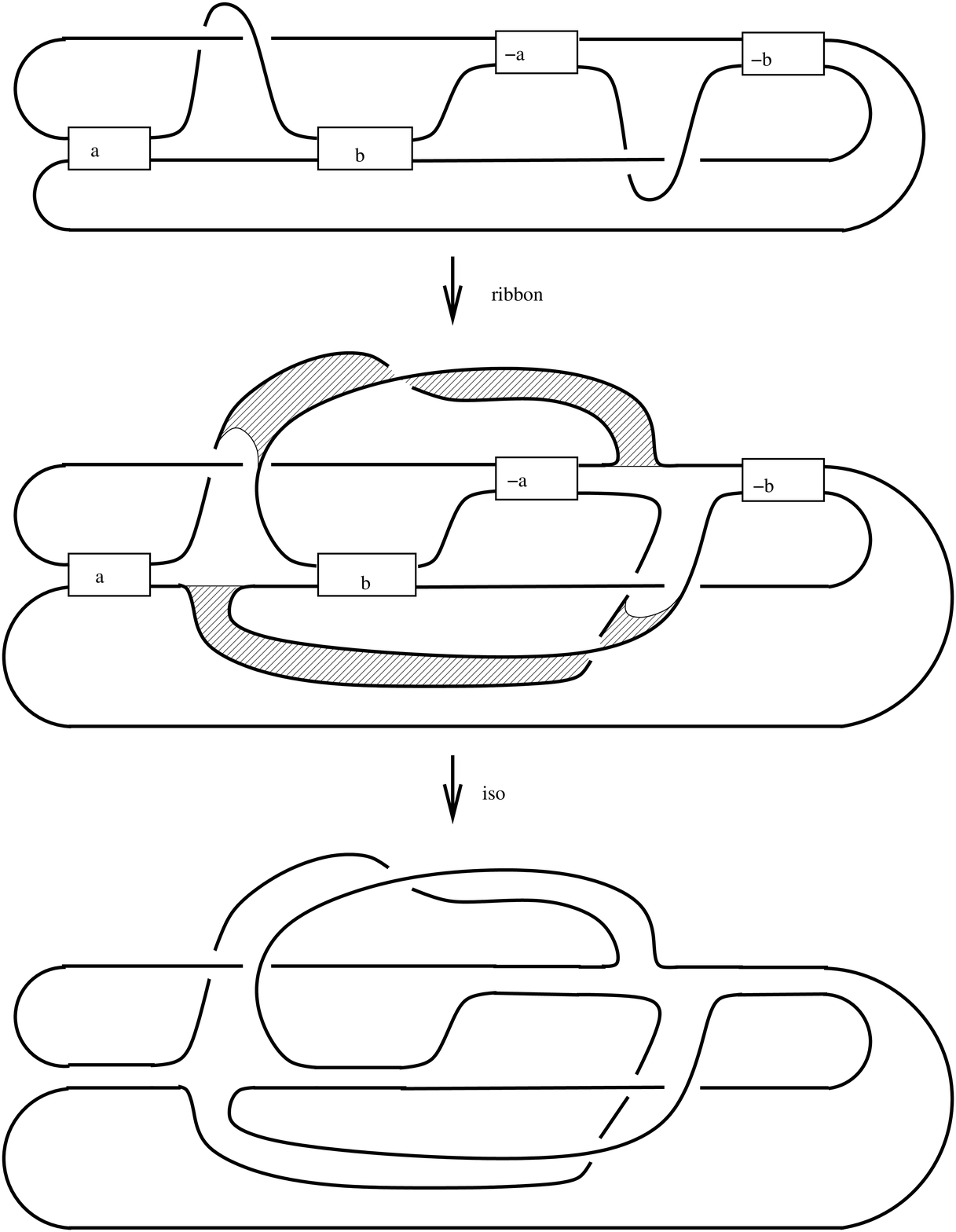}
\end{center}
\caption{Ribbon moves on the link $L_{a,b}$}
\label{f:fig3}
\end{figure}

\begin{lem}\label{l:movesI}
Let $L_{a,b}$, $a,b\in\Z$, be the 2--bridge link given by the top
picture of Figure~\ref{f:fig3}. If the link $L_{a,b}$ is a knot then
it bounds a ribbon disk. If $L_{a,b}$ has two components then it
bounds the image under a ribbon immersion of the disjoint union of a
disk and a M\"obius band.
\end{lem}

\begin{proof}
Figure~\ref{f:fig3} shows that after performing two ribbon moves the
link $L_{a,b}$ reduces to a 3--component unlink. This proves the lemma.
\end{proof}

\begin{figure}[ht]
\begin{center}
\psfrag{a}{$\scriptstyle a$}
\psfrag{b}{$\scriptstyle b$}
\psfrag{-a+1}{${\scriptstyle -a+1}$}
\psfrag{-b+1}{${\scriptstyle -b+1}$}
\psfrag{-a}{$\scriptstyle -a$}
\psfrag{b-1}{${\scriptstyle b-1}$}
\psfrag{ribbon}{(one ribbon move)}
\psfrag{iso}{(isotopy)}
\includegraphics[width=8cm]{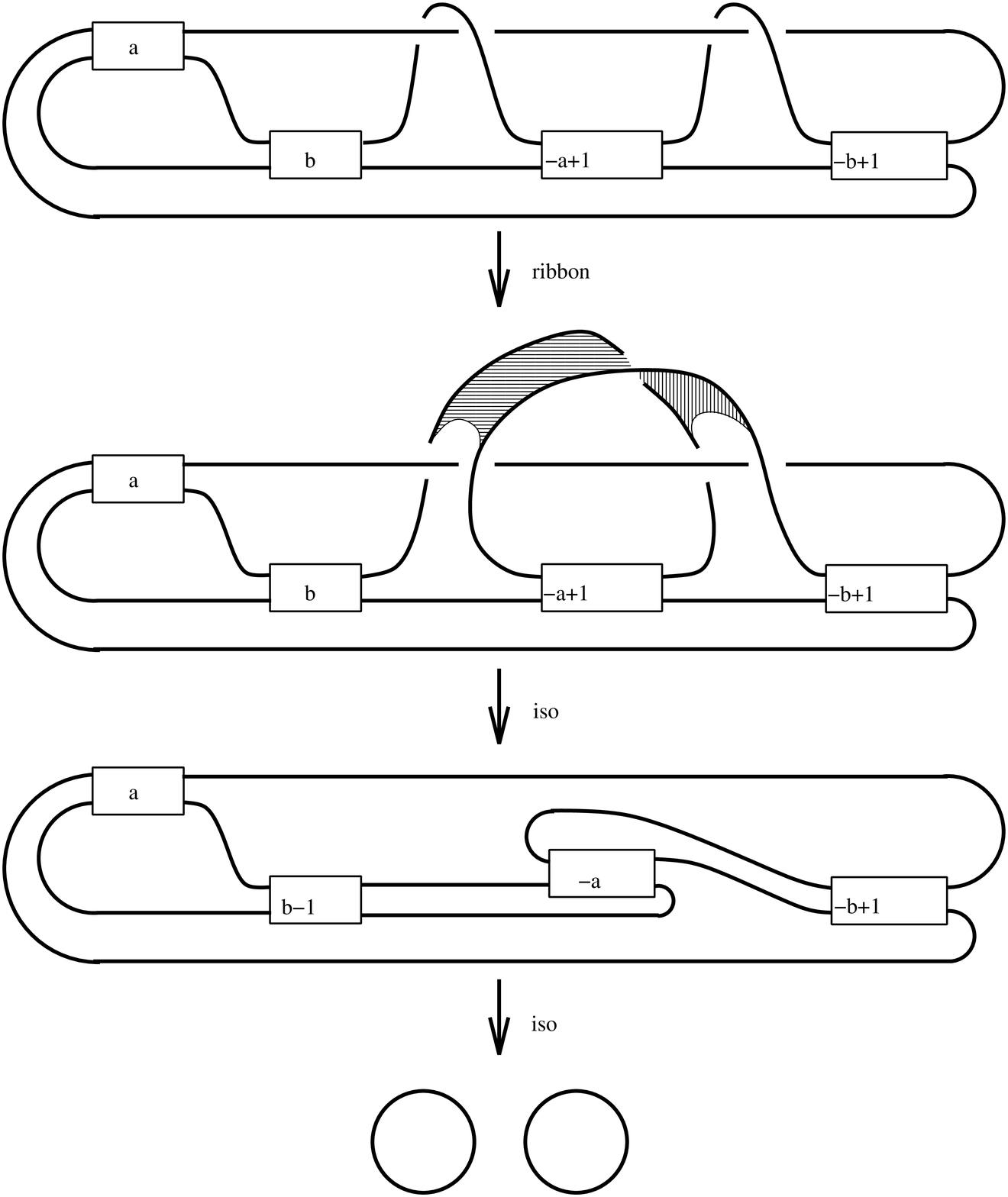}
\end{center}
\caption{Ribbon move on the link $L'_{a,b}$}
\label{f:fig4}
\end{figure}

\begin{lem}\label{l:movesII}
Let $L'_{a,b}$, $a,b\in\Z$, be the link given by the top picture of
Figure~\ref{f:fig4}. If the link $L'_{a,b}$ is a knot then it bounds a
ribbon disk. If $L'_{a,b}$ has two components then it bounds the image
under a ribbon immersion of the disjoint union of a disk and a
M\"obius band.
\end{lem}

\begin{proof}
Figure~\ref{f:fig4} shows that after performing one ribbon move the
link $L'_{a,b}$ reduces to a 2--component unlink. This proves the
lemma.
\end{proof}

\begin{lem}\label{l:familyII}
Let $p>q>0$ be coprime integers, and suppose that $\frac pq$ 
is equal to one of the following: 
\begin{enumerate}
\item
$[2^{[t]},3,s+2,t+2,3,2^{[s]}]^-$, $s,t\geq 0$,
\item
$[2^{[t]},s+3,2,t+2,3,2^{[s]}]^-$, $s,t\geq 0$.
\end{enumerate}
Then, if $p$ is odd $K(p,q)$ bounds a ribbon disk; if $p$ is even the
2--component link $K(p,q)$ bounds the image under a ribbon immersion
the disjoint union of a disk and a M\"obius band.
\end{lem}

\begin{proof}
By Equation~\eqref{e:cfrac} we have 
\[
[1,t,1,1,s,1,t,1,1,s+1]^+ = [2^{[t]},3,s+2,t+2,3,2^{[s]}]^-.
\]
Therefore, in Case $(1)$ the link $K(p,q)=L(1,t,1,1,s,1,t,1,1,s+1)$ is
given by the top picture in Figure~\ref{f:fig5}. After an isotopy, the
knot appears as in the middle picture of Figure~\ref{f:fig5}. After a
further isotopy, we obtain the bottom picture of Figure~\ref{f:fig5}.
\begin{figure}[ht]
\begin{center}
\psfrag{t}{${\scriptstyle t}$}
\psfrag{-t}{${\scriptstyle -t}$}
\psfrag{-s}{${\scriptstyle -s}$}
\psfrag{t+1}{${\scriptstyle t+1}$}
\psfrag{s+1}{${\scriptstyle s+1}$}
\psfrag{-s-1}{${\scriptstyle -s-1}$}
\psfrag{t+2}{${\scriptstyle t+2}$}
\psfrag{s+2}{${\scriptstyle s+2}$}
\psfrag{-t-1}{${\scriptstyle -t-1}$}
\psfrag{-s-1}{${\scriptstyle -s-1}$}
\psfrag{iso}{(isotopy)}
\includegraphics[width=10cm]{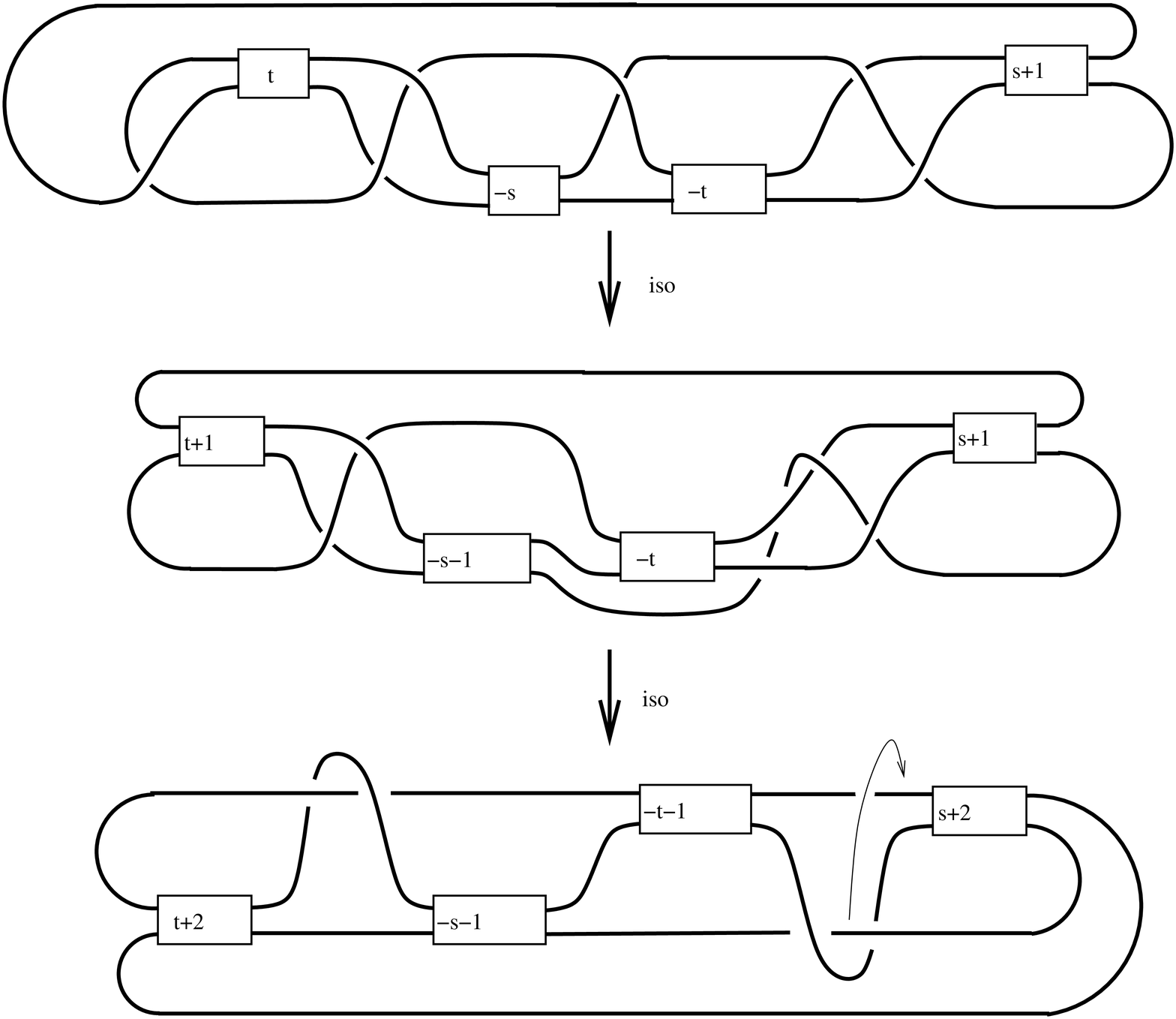}
\end{center}
\caption{The link $L(1,t,1,1,s,1,t,1,1,s+1)$}
\label{f:fig5}
\end{figure}
After an isotopy starting with pulling a strand as suggested by the
arrow in Figure~\ref{f:fig5}, we arrive at the link given by the top
picture of Figure~\ref{f:fig3} for $(a,b)=(t+2,-s-1)$. Thus, in Case
(1) the lemma follows from Lemma~\ref{l:movesI}.

By Equation~\eqref{e:cfrac} we have 
\[
[1,t,s+1,2,t,1,1,s+1]^+ = [2^{[t]},s+3,2,t+2,3,2^{[s]}]^-.
\]
Therefore $K(p,q)$ is isotopic to the link $L(1,t,s+1,2,t,1,1,s+1)$
given in Figure~\ref{f:fig6}. Applying an obvious isotopy it is easy
to see that this link is isotopic to the link $L'_{t+1,-s-1}$, where
$L'_{a,b}$, for $a,b\in\Z$, is as in Lemma~\ref{l:movesII}. Part (2)
of the statement now follows from Lemma~\ref{l:movesII}.
\end{proof} 
\begin{figure}[ht]
\begin{center}
\psfrag{t}{${\scriptstyle t}$}
\psfrag{s+1}{${\scriptstyle s+1}$}
\psfrag{-t}{${\scriptstyle -t}$}
\psfrag{-s-1}{${\scriptstyle -s-1}$}
\includegraphics[width=9cm]{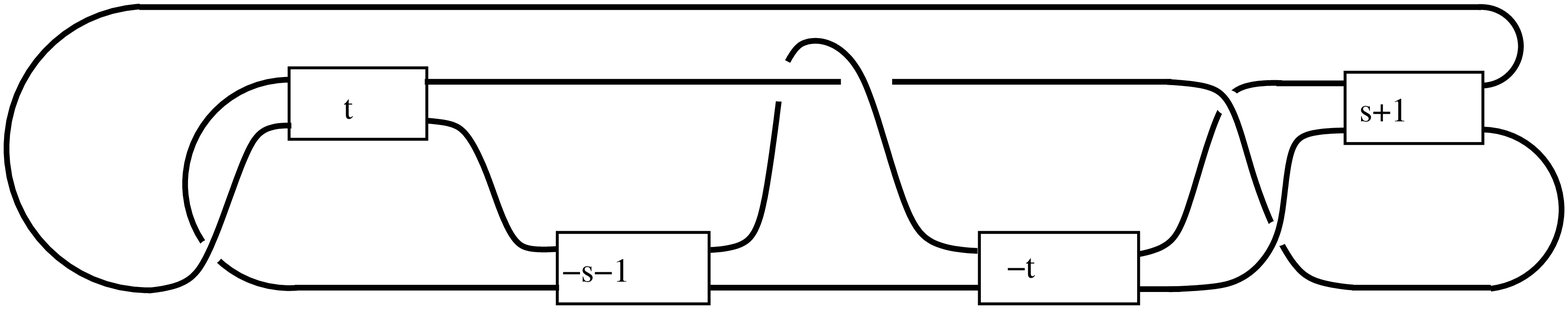}
\end{center}
\caption{The link $L(1,t,s+1,2,t,1,1,s+1)$}
\label{f:fig6}
\end{figure}

\begin{lem}\label{l:familyIII}
Let $p>q\geq 1$ be coprime integers, and suppose that $\frac pq$ 
is equal to one of the following: 
\begin{enumerate}
\item
$[t+2,s+2,3,2^{[t]},4,2^{[s]}]^-$, $s,t\geq 0$, 
\item
$[t+2,2,s+3,2^{[t]},4,2^{[s]}]^-$, $s,t\geq 0$,
\item
$[t+3,2,s+3,3,2^{[t]},3,2^{[s]}]^-$, $s,t\geq 0$.
\end{enumerate}
Then, if $p$ is odd $K(p,q)$ bounds a ribbon disk; if $p$ is even the
2--component link $K(p,q)$ bounds the image under a ribbon
immersion of the disjoint union of a disk and a M\"obius band.
\end{lem}

\begin{proof}
By Equation~\eqref{e:cfrac} we have 
\[
[t+1,1,s,1,1,t+1,2,s+1]^+ = 
[t+2,s+2,3,2^{[t]},4,2^{[s]}]^-,\quad s,t\geq 0.
\]
Therefore in Case (1) $K(p,q)$ is isotopic to the link
$L(t+1,1,s,1,1,t+1,2,s+1)$ given by the top picture of
Figure~\ref{f:fig7}.
\begin{figure}[ht]
\begin{center}
\psfrag{-t-1}{${\scriptstyle -t-1}$}
\psfrag{s+1}{${\scriptstyle s+1}$}
\psfrag{t+1}{${\scriptstyle t+1}$}
\psfrag{-s}{${\scriptstyle -s}$}
\psfrag{-t-2}{${\scriptstyle -t-2}$}
\psfrag{-s-2}{${\scriptstyle -s-2}$}
\psfrag{iso}{(isotopy)}
\includegraphics[width=9cm]{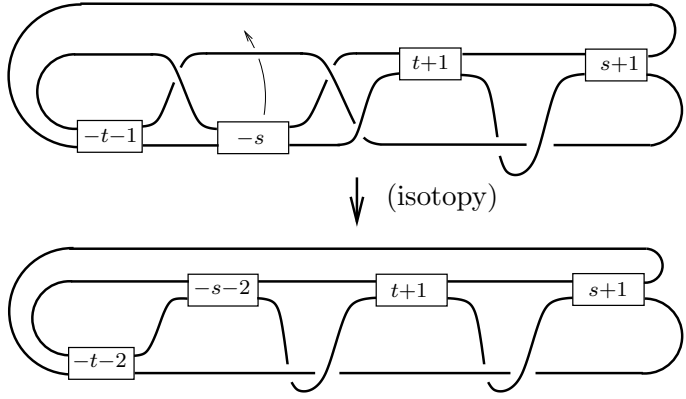}
\end{center}
\caption{The link $L(t+1,1,s,1,1,t+1,2,s+1)$}
\label{f:fig7}
\end{figure}
Applying the isotopy suggested by the arrow one obtains the link given
by the bottom picture of Figure~\ref{f:fig7}, which is easily checked
to be the mirror image of the link $L'_{t+2,s+2}$, where $L'_{a,b}$, for
$a,b\in\Z$, is as in Lemma~\ref{l:movesII}. Therefore, Part (1) of the
statement follows from Lemma~\ref{l:movesII}. By
Equation~\eqref{e:cfrac} we have
\[
[t+1,2,s+1,t+1,2,s+1]^+ = 
[t+2,2,3+s,2^{[t]},4,2^{[s]}]^-,\quad s,t\geq 0.
\]
This shows that in Case (2) $K(p,q)$ is isotopic to
$L(t+1,2,s+1,t+1,2,s+1)$, which is easily seen to be isotopic
to $L_{-t-1,-s-1}$, where $L_{a,b}$, $a,b\in\Z$ is as in
Lemma~\ref{l:movesI}. Thus, in Case (2) the statement follows from
Lemma~\ref{l:movesI}. Now observe that if
\[
\frac pq = [t+3,2,s+3,3,2^{[t]},3,2^{[s]}]^-,\quad s,t\geq 0,
\]
then by Equations~\eqref{e:Riemen} and~\eqref{e:inversion} we have
\[
\frac p{p-q'} = [s+2,t+3,3,2^{[s]},4,2^{[t+1]}]^-,
\]
which is of the type considered in Case (1). This concludes the proof. 
\end{proof}

\section{The proof of Theorem~\ref{t:main}}
\label{s:final}

\textbf{Outline.}  In this section we use the results obtained in the
previous sections to prove Theorem~\ref{t:main}. 

Before starting the proof of Theorem~\ref{t:main} we need four
arithmetic lemmas.

\begin{lem}\label{l:fraction}
Suppose that $a_i\geq 2$ for $i=1,\ldots, n$, are integers and
\[
[a_1,\ldots,a_n]^- = \frac{m^2}{mk\pm 1},\quad (m,k)=1,\quad  0< k<m.
\]
Then,
\[
[2,a_1,\ldots,a_n,a_n + 1]^- = \frac{(2m-k)^2}{(2m-k)m\pm 1} 
\]
and
\[
[a_1+1,a_2,\ldots,a_n,2]^- = \frac{(m+k)^2}{(m+k)k\pm 1}.
\]
\end{lem}

\begin{proof}
Since
\[
(m(m-k)\pm 1)(mk\pm 1) \equiv 1\pmod {m^2},
\]
by Equation~\eqref{e:inversion} we have
\begin{equation}\label{e:reverse}
[a_n,\ldots,a_1]^-=\frac{m^2}{m(m-k)\pm 1},
\end{equation}
therefore 
\[
[a_n + 1,\ldots,a_1]^-=\frac{m^2+m(m-k)\pm 1}{m(m-k)\pm 1}
=\frac{2m^2-mk\pm 1}{m(m-k)\pm 1}.
\]
Similarly, since
\[
(m(m-k)\pm 1)(2mk-k^2\pm 2) = 1 - (2m^2-mk\pm 1)(k^2-mk\mp 1)
\]
we have 
\[
[a_1,\ldots,a_n +1]^-=\frac{2m^2-mk\pm 1}{2mk-k^2\pm 2}.
\]
The first formula in the statement of the lemma now follows by a simple
computation. By Equation~\eqref{e:reverse} and the first formula 
in the statement we have 
\[
[2,a_n,\ldots,a_1+1]^-=\frac{(2m-(m-k))^2}{(2m-m+k)m\pm 1} = 
\frac{(m+k)^2}{(m+k)m\pm 1},
\]
which implies, as before, the second formula in the statement of the lemma.
\end{proof}

\begin{lem}\label{l:I=-3(3)}
Let $n\geq 3$ and let $S_n=\{v_1,\ldots,v_n\}\subseteq\bD^n$ be a
standard subset such that $I(S_n)=-3$. Suppose $v_i\cdot v_i=-a_i$ for
$i=1,\ldots,n$. Then, 
\[
[a_1,\ldots,a_n]^-=
\frac{m^2}{mk + 1}
\]
for some integers $m,k$ with $0<k<m$ and $(m,k)=1$.
\end{lem}

\begin{proof}
The fraction associated to the set $S_3$ of Lemma~\ref{l:n=3} is
$[2,2,2]=4/3$, which is of the form $m^2/(m+1)$. The lemma follows
immediately from Lemmas~\ref{l:I=-3(2)} and~\ref{l:fraction}.
\end{proof}

In the following proofs we shall use the formula 
\begin{equation}\label{e:utile}
[2^{[t]},x]^-=\frac{(t+1)x-t}{tx-(t-1)},\quad
t\in\N\cup\{0\},
\end{equation}
which holds for any variable $x$ and can be established by an easy
induction.

\begin{lem}\label{l:I=-2(3)}
Let $n\geq 4$, and let $S_n=\{v_1,\ldots,v_n\}\subseteq\bD^n$ be a
standard subset such that $I(S_n)=-2$. Suppose $v_i\cdot v_i=-a_i$ for
$i=1,\ldots,n$. Then, either $[a_1,\ldots,a_n]^-$ or
$[a_n,\ldots,a_1]^-$ is of one of the following forms:
\begin{enumerate}
\item
$\dfrac{m^2}{m^2-d(m-1)}$, where $d$ divides $2m+1$; 
\item
$\dfrac{m^2}{m^2-d(m-1)}$, where $d$ is odd and divides $m-1$.
\end{enumerate}
\end{lem}

\begin{proof}
Using Formula~\eqref{e:utile} one can verify
that
\begin{multline*}
[2^{[t]},3,s+2,t+2,3,2^{[s]}]^- =\\ 
\frac{(2st+3s+3t+4)^2}{(2st+3s+3t+4)^2-(2s+3)(2st+3s+3t+3)}.
\end{multline*}
For $m=2st+3s+3t+4$ and $d=2s+3$, since $2m+1=(2s+3)(2t+3)$, 
this shows that the associated fraction is of the form 
\[
\frac{m^2}{m^2-d(m-1)},
\]
where $d$ divides $2m+1$. Similarly, 
\begin{multline*}
[2^{[t]},s+3,2,t+2,3,2^{[s]}]^- = 
\frac{(2st+2s+3t+4)^2}{(2st+2s+3t+4)^2-(2s+3)^2 (t+1)}=\\
\frac{m^2}{m^2-d(m-1)},
\end{multline*}
where $m=2st+2s+3t+4$ and $d=2s+3$. Observe that $d$ is odd and divides
$m-1=(2s+3)(t+1)$. By Lemma~\ref{l:I=-2(2)} this concludes the proof.
\end{proof}

\begin{lem}\label{l:I=-1(2)}
Let $n\geq 4$ and let $S_n=\{v_1,\ldots,v_n\}\subseteq\bD^n$ be a
standard subset such that $I(S_n)=-1$. Suppose that $v_i\cdot
v_i=-a_i$, $i=1,\ldots,n$. Then, either $[a_1,\ldots,a_n]^-$ or
$[a_n,\ldots,a_1]^-$ is of one of the following types:
\begin{enumerate}
\item
$\dfrac{m^2}{d(m+1)}$ where $d$ is odd and divides $m+1$; 
\item
$\dfrac{m^2}{d(m+1)}$ where $d$ divides $2m-1$;
\item
$\dfrac{m^2}{m^2-d(m+1)}$ where $d$ is odd and divides $m+1$.
\end{enumerate}
\end{lem}

\begin{proof}
Using Formula~\eqref{e:utile} one can verify that
\[
[t+2,s+2,3,2^{[t]},4,2^{[s]}]^- = 
\frac{(2st+4s+3t+5)^2}{(2s+3)^2(t+2)} = 
\frac{m^2}{d(m+1)},
\]
where $m=2st+4s+3t+5$ and $d=2s+3$ is odd and divides
$m+1=(2s+3)(t+2)$. Similarly,
\[
[t+2,2,s+3,2^{[t]},4,2^{[s]}]^- = 
\frac{(2st+3s+3t+5)^2}{(2s+3)(2st+3s+3t+6)}=
\frac{m^2}{d(m+1)},
\]
where $m=2st+3s+3t+5$ and $d=2s+3$ divides $2m-1=(2s+3)(2t+3)$. 
Finally, 
\begin{multline*}
[t+3,2,s+3,3,2^{[t]},3,2^{[s]}]^- =
\frac{(2ts+5s+4t+9)^2}{(s+2)(4ts+10s+8t+17)}= \\
\frac{m^2}{(2m-1)(m+1)/d},
\end{multline*}
where $m=2ts+5s+4t+9$ and $d=2t+5$ divides $m+1=(s+2)(2t+5)$. Since 
\[
\frac{(2m-1)(m+1)}d (m^2-d(m+1)) \equiv 1\pmod {m^2},
\]
by Equation~\eqref{e:inversion} we have 
\[
\frac{m^2}{m^2-d(m+1)} = [2^{[s]},3,2^{[t]},3,s+3,2,t+3]^-.
\]
Thus, the lemma follows by Lemma~\ref{l:I=-1}.
\end{proof}

\begin{proof}[Proof of Theorem~\ref{t:main}]
We first show that (2) implies (1). Let us assume for that (2)
holds. Let $\widetilde\Si\subset B^4$ be a smoothly embedded surface
obtained by pushing the interior of $\Si$ inside the 4--ball. It is
easy to check that (regardless of the parity of $p$) the inclusion
$S^3\setminus\del\widetilde\Si\subset B^4\setminus\widetilde\Si$
induces a surjective homomorphism
\[
\varphi\co
H_1(S^3\setminus\del\widetilde\Si;\Z)\to H_1(B^4\setminus\widetilde\Si;\Z)
\]
such that the homomorphism
$H_1(S^3\setminus\del\widetilde\Si;\Z)\to\Z/2\Z$ defining the 2--fold
cover $L(p,q)\to S^3$ branched along $\del\widetilde\Si=K(p,q)$
factors through $H_1(B^4\setminus\widetilde \Si;\Z)$ via
$\varphi$. Therefore, the cover $L(p,q)\to S^3$ extends to a 2--fold
cover $W\to B^4$ branched along $\widetilde\Si$. We may assume that
the distance function from the origin $B^4\to [0,1]$ restricted to
$\widetilde\Si$ is a proper Morse function with only index--0 and
index--1 critical points. This implies that $W$ has a handlebody
decomposition with only 0--, 1-- and 2--handles (see e.g.~\cite[lemma
at pages 30--31]{CH}). Therefore, from
\[
b_0(W)-b_1(W)+b_2(W) = \chi(W) = 2\chi(B^4) - \chi(\widetilde\Si) = 1
\]
we deduce $b_1(W)=b_2(W)$. On the other hand, since $b_1(\del W)=0$
and $H_1(W,\del W;\Q)\cong H^3(W;\Q)=0$, the homology exact sequence
of the pair $(W,\del W)$ gives $b_1(W)=0$, so it follows that
$H_*(W;\Q)\cong H_*(B^4;\Q)$, and (1) holds. 

Now we show that (1) implies (3). Assume that Part (1) of the
statement holds. It is a well--known fact that if $\frac pq =
[a_1,\ldots,a_n]^-$ the lens space $L(p,q)$ smoothly bounds the
4--dimensional plumbing $P(p,q)$ given by the weighted graph of
Figure~\ref{f:fig8}.
\begin{figure}[ht]
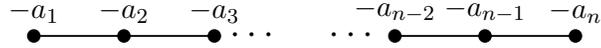

\begin{center}
\setlength{\unitlength}{1mm}
\unitlength=1.2cm
\begin{graph}(8,0.5)(0,0)
\graphnodesize{0.15}
  \roundnode{n1}(1,0)
  \roundnode{n2}(2,0)
  \roundnode{n3}(3,0)
  \roundnode{n4}(5,0)
  \roundnode{n5}(6,0)
  \roundnode{n6}(7,0)

  \edge{n1}{n2}
  \edge{n2}{n3}
  \edge{n4}{n5}
  \edge{n5}{n6}

  \autonodetext{n1}[n]{$-a_1$}
  \autonodetext{n2}[n]{$-a_2$}
  \autonodetext{n3}[n]{$-a_3$}
  \autonodetext{n3}[e]{{\Large$\cdots$}}
  \autonodetext{n4}[w]{{\Large$\cdots$}}
  \autonodetext{n4}[n]{$-a_{n-2}$}
  \autonodetext{n5}[n]{$-a_{n-1}$}
  \autonodetext{n6}[n]{$-a_n$}
\end{graph}
\end{center}
\caption{The graph of the canonical plumbing bounded by $L(p,q)$}
\label{f:fig8}
\end{figure}
The intersection form of $P(p,q)$ is negative definite. Hence, since
$L(p,q)\cong -L(p,p-q)$, if $L(p,q)$ smoothly bounds a rational
homology 4--ball $W(p,q)$ we can construct the smooth, negative
4--manifolds
\[
X(p,q)=P(p,q)\cup_{\del} (-W(p,q)),\quad 
X(p,p-q)=P(p,p-q)\cup_{\del} W(p,q)
\]
By Donaldson's theorem on the intersection form of definite
4--manifolds~\cite{Do}, the intersection forms of $X(p,q)$ and
$X(p,p-q)$ are both standard diagonal. Hence, suppose that the
intersection lattice of $X(p,q)$ is isomorphic to $\bD^n$ and the
intersection lattice of $X(p,p-q)$ is isomorphic to
$\bD^{n'}$. Clearly, the intersection lattices
$H_2(P(p,q);\Z)\cong\Z^n$ and $H_2(P(p,p-q);\Z)\cong\Z^{n'}$ have
bases $\{v_1,\ldots,v_n\}$ and $\{w_1,\ldots,w_{n'}\}$ which satisfy
Equations~\eqref{e:inters0}. Therefore, via the embeddings
$P(p,q)\subset X(p,q)$ and $P(p,p-q)\subset X(p,p-q)$ we can view the
above bases as standard subsets $S\subset\bD^n$ and
$S'\subset\bD^{n'}$ with associated strings $(a_1,\ldots,a_n)$ and
$(b_1,\ldots,b_{n'})$, where $[b_1,\ldots,b_{n'}]^-=p/(p-q)$. In view
of Lemma~\ref{l:negsum}, we may assume without loss of generality that
$I(S)<0$. Then, by Theorem~\ref{t:standard} and Lemmas~\ref{l:n=3},
\ref{l:I=-3(3)}, \ref{l:I=-2(3)} and \ref{l:I=-1(2)} it follows that
(3) holds.

Finally, we show that (3) implies (2). Suppose that (3) holds,
i.e.~$\frac pq\in\RR$. Then, since applying finitely many times the
functions $f$ and $g$ of Definition~\ref{d:R} amounts to changing
$K(p,q)$ by an isotopy or a reflection, we may assume that $p=m^2$ and
$q$ is of one of the three types given in Definition~\ref{d:R}. We
consider various cases separately.

{\em First case: $q=mk\pm 1$, with $m>k>0$ and $(m,k)=1$.}

In view of Lemmas~\ref{l:I=-3(2)} and~\ref{l:familyI}, it suffices to
show that the string of coefficients of the continued fraction expansion
of $\frac pq$ is obtained from $(2,2,2)$ via a finite sequence of
operations as in Lemma~\ref{l:I=-3(2)}. Since $m^2-(mk\mp 1)=m(m-k)\pm 1$
and either $m\geq 2k$ or $m\geq 2(m-k)$, up to replacing $k$ with
$m-k$ (and $K(p,q)$ with its mirror image $K(p,p-q)$) we may assume
$m\geq 2k$. If $m=2k$, since $(m,k)=1$ we must have $m=2$, $k=1$ and
$p/q=[2,2,2]^-$. If $m>2k$, arguing by induction on $m$ we may assume
\[
\frac{(m-k)^2}{(m-k)k\pm 1} = [a_1,a_2,\ldots,a_n]^-,
\]
where $(a_1,a_2,\ldots,a_n)$ is obtained from $(2,2,2)$ as described
above. But in view of Lemma~\ref{l:fraction} we have
\[
\frac{m^2}{mk\pm 1} = [a_1+1,a_2,\ldots,a_n,2]^-,
\]
so we are done. 

{\em Second case: $q=d(m - 1)$, where $d>1$ divides $2m + 1$.} 

It suffices to show that (2) holds for $K(p,p-q)$. Since $d(m-1)<m^2$,
we have $2m+1>d>1$, and $d$ must be odd because it divides $2m+1$.
Therefore we can write $d=2s+3$ for some $s\geq 0$ and $2m+1=d(2t+3)$
for some $t\geq 0$. Then $m=2st+3s+3t+4$, and as in the proof of 
Lemma~\ref{l:I=-2(3)}
\[
\frac{m^2}{m^2-d(m-1)} = [2^{[t]},3,s+2,t+2,3,2^{[s]}]^-.
\]
Therefore (2) holds by Lemma~\ref{l:familyII}(1). 

{\em Third case: $q=d(m+1)$, where $d>1$ divides $2m-1$.}

Arguing as in the previous case, we can write $d=2s+3$ and
$2m-1=d(2t+3)$ for some $s,t\geq 0$. Then, $m=2st+3s+3t+5$ and 
\[
\frac{m^2}{d(m+1)} = [t+2,2,s+3,2^{[t]},4,2^{[s]}]^-,
\]
which implies (2) by Lemma~\ref{l:familyIII}(2). 

{\em Fourth case: $q=d(m+1)$, where $d>1$ is odd and divides $m+1$.}

Since $d(m+1)<m^2$ we have $m+1>d>1$, therefore we can write $d=2s+3$ 
and $m+1=d(t+2)$ for some $s,t\geq 0$. Then
\[
\frac{m^2}{d(m+1)} = [t+2,s+2,3,2^{[t]},4,2^{[s]}]^-,
\]
and (2) holds by Lemma~\ref{l:familyIII}(1).

{\em Fifth case: $q=d(m-1)$, where $d>1$ is odd and divides $m-1$.}

As before, it suffices to prove that (2) holds for $K(p,p-q)$. We can
write $d=2s+3$ and $m-1=d(t+1)$ for some $s,t\geq 0$.  Then
\[
\frac{m^2}{m^2 - d(m-1)} = [2^{[t]},s+3,2,t+2,3,2^{[s]}]^-,
\]
and (2) holds by Lemma~\ref{l:familyII}(2). This concludes the proof.
\end{proof}

\end{document}